\documentclass[reqno,12pt]{amsart}
\theoremstyle{plain}    
\newtheorem{thm}{Theorem}[section]
\theoremstyle{plain}    
\newtheorem{lemma}[thm]{Lemma} 
\newtheorem{prop}[thm]{Proposition}
\theoremstyle{remark}
\newtheorem{defi}[thm]{Definition}

\newtheorem{example}[thm]{Example}
\newtheorem{examples}[thm]{Examples}
\newtheorem{notation}[thm]{Notation}
\newtheorem{remark}[thm]{Remark}

\newtheorem{ques}[thm]{Question}

\usepackage{amscd,amssymb,comment,epic,euscript,graphics}

\newcommand\Afr{{\mathfrak A}}

\newcommand\alg{{\operatorname{alg}}}

\newcommand\alphat{{\tilde\alpha}}

\newcommand\Bc{{\mathcal{B}}}

\newcommand\betat{{\tilde\beta}}

\newcommand\Cpx{{\mathbf C}}

\newcommand\Dc{{\mathcal{D}}}

\newcommand\diag{\text{\rm diag}}

\newcommand\DT{\operatorname{DT}}

\newcommand\Ec{{\mathcal{E}}}

\newcommand\eps{\epsilon}

\newcommand\Fb{{\mathbf F}}

\newcommand\Gt{{\widetilde G}}

\newcommand\HEu{{\EuScript H}}                   

\newcommand\Ints{{\mathbf Z}}

\newcommand\KEu{{\EuScript K}}                   

\newcommand\Mcal{{\mathcal{M}}} 
\newcommand\MEu{{\EuScript M}}                   

\newcommand\Nats{{\mathbf N}}

\newcommand\Nc{{\mathcal{N}}}

\newcommand\NC{\operatorname{NC}}

\newcommand\oneh{{\hat 1}}

\newcommand\phit{{\tilde\phi}}

\newcommand\pit{{\tilde\pi}}

\newcommand\Proj{{\mathrm{Proj}}}

\newcommand\Reals{{\mathbf R}}

\newcommand\restrict{{\upharpoonright}}

\newcommand\xh{{\hat x}}

\newcommand\Vc{{\mathcal{V}}}

\newcommand\vN{{vN}}

\newcommand\wt{{\tilde w}}

\begin{document}

\pagestyle{myheadings}

 \title{Hyperinvariant subspaces for some $B$--circular operators}

 \author{Ken Dykema}

 \address{\hskip-\parindent
 Department of Mathematics \\
 Texas A\&M University \\
 College Station TX 77843--3368, USA}
 \email{kdykema@math.tamu.edu}

 \thanks{Supported in part by a grant from the NSF}

 \date{April 19, 2004}

 \begin{abstract}
 We show that if $A$ is a Hilbert--space operator, then the set of all projections
 onto hyperinvariant subspaces of $A$, which is contained in the von Neumann algebra $\vN(A)$
 that is generated by $A$, is independent of the representation of $\vN(A)$, thought of as
 an abstract W$^*$--algebra.
 
 We modify a technique of Foias, Ko, Jung and Pearcy to get a method for finding nontrivial hyperinvariant subspaces
 of certain operators in finite von Neumann algebras.

 We introduce the $B$--circular operators as a special case of Speicher's $B$--Gaussian operators in free probability theory,
 and we prove several results about a $B$--circular operator $z$, including formulas for the $B$--valued Cauchy-- and
 R--transforms of $z^*z$.
 We show that a large class of $L^\infty([0,1])$--circular operators in finite von Neumann algebras have nontrivial
 hyperinvariant subspaces, and that another large class of them can be embedded in the free group factor
 $L(\Fb_3)$.
 These results generalize some of what is known about the quasinilpotent DT--operator.
 \end{abstract}

 \maketitle


\section{Introduction}

The invariant subspace problem for operators on Hilbert space and the related hyperinvariant subspace problem
are both unresolved and are of importance for understanding the structure of Hilbert space operators.
Let $\HEu$ be a Hilbert space and let $A\in\Bc(\HEu)$ be a bounded operator on $\HEu$.
A closed subspace $\HEu_0\subseteq\HEu$ is said to be $A$--{\em invariant} if $A(\HEu_0)\subseteq\HEu_0$.
(Throughout this paper, all subspaces will be assumed to be closed.)
The subspace $\HEu_0$ is said to be $A$--{\em hyperinvariant}
if it is $S$--invariant whenever $S\in\Bc(\HEu)$ commutes with $A$.
Recall that the invariant subspace problem asks whether, for $\HEu$ infinite dimensional, every $A\in\Bc(\HEu)$
has an $A$--invariant subspace that is nontrivial (i.e.\ neither $\{0\}$ nor $\HEu$ itself), and
the hyperinvariant subspace problem asks whether every $A\in\Bc(\HEu)$ that is not a scalar multiple of the identity
has a nontrivial $A$--hyperinvariant subspace.

\begin{figure}[t]
\begin{picture}(80,80)(0,0)

  \linethickness{0.7pt}
  \drawline(0,0)(0,80)(80,80)(80,0)(0,0)

  \linethickness{0.5pt}
  \drawline(0,80)(80,0)

  \linethickness{0.3pt}

  \drawline(2,78)(80,78)
  \drawline(4,76)(80,76)
  \drawline(6,74)(80,74)
  \drawline(8,72)(80,72)
  \drawline(10,70)(80,70)
  \drawline(12,68)(80,68)
  \drawline(14,66)(80,66)
  \drawline(16,64)(80,64)
  \drawline(18,62)(80,62)
  \drawline(20,60)(80,60)
  \drawline(22,58)(80,58)
  \drawline(24,56)(80,56)
  \drawline(26,54)(80,54)
  \drawline(28,52)(80,52)
  \drawline(30,50)(80,50)
  \drawline(32,48)(80,48)
  \drawline(34,46)(80,46)
  \drawline(36,44)(80,44)
  \drawline(38,42)(80,42)
  \drawline(40,40)(80,40)
  \drawline(42,38)(80,38)
  \drawline(44,36)(80,36)
  \drawline(46,34)(80,34)
  \drawline(48,32)(80,32)
  \drawline(50,30)(80,30)
  \drawline(52,28)(80,28)
  \drawline(54,26)(80,26)
  \drawline(56,24)(80,24)
  \drawline(58,22)(80,22)
  \drawline(60,20)(80,20)
  \drawline(62,18)(80,18)
  \drawline(64,16)(80,16)
  \drawline(66,14)(80,14)
  \drawline(68,12)(80,12)
  \drawline(70,10)(80,10)
  \drawline(72,8)(80,8)
  \drawline(74,6)(80,6)
  \drawline(76,4)(80,4)
  \drawline(78,2)(80,2)

  \drawline(2,78)(2,80)
  \drawline(4,76)(4,80)
  \drawline(6,74)(6,80)
  \drawline(8,72)(8,80)
  \drawline(10,70)(10,80)
  \drawline(12,68)(12,80)
  \drawline(14,66)(14,80)
  \drawline(16,64)(16,80)
  \drawline(18,62)(18,80)
  \drawline(20,60)(20,80)
  \drawline(22,58)(22,80)
  \drawline(24,56)(24,80)
  \drawline(26,54)(26,80)
  \drawline(28,52)(28,80)
  \drawline(30,50)(30,80)
  \drawline(32,48)(32,80)
  \drawline(34,46)(34,80)
  \drawline(36,44)(36,80)
  \drawline(38,42)(38,80)
  \drawline(40,40)(40,80)
  \drawline(42,38)(42,80)
  \drawline(44,36)(44,80)
  \drawline(46,34)(46,80)
  \drawline(48,32)(48,80)
  \drawline(50,30)(50,80)
  \drawline(52,28)(52,80)
  \drawline(54,26)(54,80)
  \drawline(56,24)(56,80)
  \drawline(58,22)(58,80)
  \drawline(60,20)(60,80)
  \drawline(62,18)(62,80)
  \drawline(64,16)(64,80)
  \drawline(66,14)(66,80)
  \drawline(68,12)(68,80)
  \drawline(70,10)(70,80)
  \drawline(72,8)(72,80)
  \drawline(74,6)(74,80)
  \drawline(76,4)(76,80)
  \drawline(78,2)(78,80)

\end{picture}
\caption{The quasinilpotent DT--operator $T$} \label{fig:T}
\end{figure}
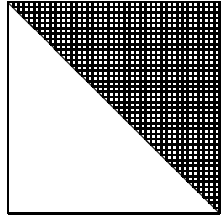

Uffe Haagerup~\cite{Ha} made a huge advance on the hyperinvariant subspace problem
for operators in II$_1$--factors.
He proved that if $A$ belongs to a II$_1$--factor that is embeddable in 
the ultrapower $R^\omega$ of the hyperfinite II$_1$--factor
and if the Brown measure~\cite{Br} of $A$ is supported on more than one point,
then $A$ has a nontrivial hyperinvariant subspace.
(He actually proved much more, namely a  result on Brown measure decomposition by restricting to hyperinvariant subspaces.)
It is therefore of particular interest to study the hyperinvariant subspace problem for operators
whose Brown measure has support reduced to a single point.
Since the support of the Brown measure is contained in the spectrum of the operator, quasinilpotent operators
in II$_1$--factors are of special interest.
The quasinilpotent DT--operator $T$ in the free group factor $L(\Fb_2)$,
from the family of operators defined in~\cite{DH2}, was a particularly compelling
example to study.
The operator $T$ can be realized as a limit in $*$--moments of strictly upper triangular random matrices
with i.i.d.\ complex Gaussian entries above the diagonal.
Alternatively, as was seen in~\cite[\S4]{DH2}, $T$ can be obtained from a semicircular element $X$
and a free copy of $L^\infty([0,1])$ by using projections from the latter to cut out the upper triangular 
part of $X$; for future reference, note that, $X$ may be replaced by a circular operator for this procedure.
Pictorially, then, we may represent $T$ as in Figure~\ref{fig:T}.
Here the shaded region has weight $1$, the unshaded region has weight $0$, and these weights are used to multiply
entries of a Gaussian random matrix, as was similarly considered in the
self--adjoint case by Shlyakhtenko in~\cite{Sh96} and~\cite{Sh98}.

In~\cite{DH3}, Haagerup and the author proved that $T$ has a one--parameter family of nontrivial hyperinvariant
subspaces.
The proof utilized precise knowledge of certain $*$--moments of $T$, conjectured
in~\cite{DH2} and proved by \'Sniady~\cite{Sn}, which implies that $TT^*$ and $k(T^k(T^*)^k)^{1/k}$
have the same moments for every $k\in\Nats$.
It was also shown in~\cite{DH3} that these hyperinvariant subspaces can be characterized in
terms of the asymptotic rate of decay of $\|T^n\xi\|$ as $n\to\infty$, for vectors $\xi$ in the Hilbert space.

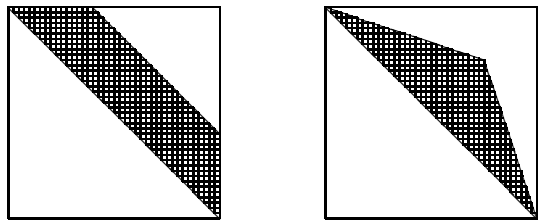
\begin{figure}[b]
\centering
\begin{picture}(200,90)

  \linethickness{0.7pt}
  \drawline(0,0)(0,80)(80,80)(80,0)(0,0)

  \linethickness{0.5pt}
  \drawline(0,80)(80,0)
  \drawline(32,80)(80,32)

  \linethickness{0.3pt}

  \drawline(2,78)(34,78)
  \drawline(4,76)(36,76)
  \drawline(6,74)(38,74)
  \drawline(8,72)(40,72)
  \drawline(10,70)(42,70)
  \drawline(12,68)(44,68)
  \drawline(14,66)(46,66)
  \drawline(16,64)(48,64)
  \drawline(18,62)(50,62)
  \drawline(20,60)(52,60)
  \drawline(22,58)(54,58)
  \drawline(24,56)(56,56)
  \drawline(26,54)(58,54)
  \drawline(28,52)(60,52)
  \drawline(30,50)(62,50)
  \drawline(32,48)(64,48)
  \drawline(34,46)(66,46)
  \drawline(36,44)(68,44)
  \drawline(38,42)(70,42)
  \drawline(40,40)(72,40)
  \drawline(42,38)(74,38)
  \drawline(44,36)(76,36)
  \drawline(46,34)(78,34)
  \drawline(48,32)(80,32)
  \drawline(50,30)(80,30)
  \drawline(52,28)(80,28)
  \drawline(54,26)(80,26)
  \drawline(56,24)(80,24)
  \drawline(58,22)(80,22)
  \drawline(60,20)(80,20)
  \drawline(62,18)(80,18)
  \drawline(64,16)(80,16)
  \drawline(66,14)(80,14)
  \drawline(68,12)(80,12)
  \drawline(70,10)(80,10)
  \drawline(72,8)(80,8)
  \drawline(74,6)(80,6)
  \drawline(76,4)(80,4)
  \drawline(78,2)(80,2)

  \drawline(2,78)(2,80)
  \drawline(4,76)(4,80)
  \drawline(6,74)(6,80)
  \drawline(8,72)(8,80)
  \drawline(10,70)(10,80)
  \drawline(12,68)(12,80)
  \drawline(14,66)(14,80)
  \drawline(16,64)(16,80)
  \drawline(18,62)(18,80)
  \drawline(20,60)(20,80)
  \drawline(22,58)(22,80)
  \drawline(24,56)(24,80)
  \drawline(26,54)(26,80)
  \drawline(28,52)(28,80)
  \drawline(30,50)(30,80)
  \drawline(32,48)(32,80)
  \drawline(34,46)(34,78)
  \drawline(36,44)(36,76)
  \drawline(38,42)(38,74)
  \drawline(40,40)(40,72)
  \drawline(42,38)(42,70)
  \drawline(44,36)(44,68)
  \drawline(46,34)(46,66)
  \drawline(48,32)(48,64)
  \drawline(50,30)(50,62)
  \drawline(52,28)(52,60)
  \drawline(54,26)(54,58)
  \drawline(56,24)(56,56)
  \drawline(58,22)(58,54)
  \drawline(60,20)(60,52)
  \drawline(62,18)(62,50)
  \drawline(64,16)(64,48)
  \drawline(66,14)(66,46)
  \drawline(68,12)(68,44)
  \drawline(70,10)(70,42)
  \drawline(72,8)(72,40)
  \drawline(74,6)(74,38)
  \drawline(76,4)(76,36)
  \drawline(78,2)(78,34)

  \linethickness{0.7pt}
  \drawline(120,0)(120,80)(200,80)(200,0)(120,0)

  \linethickness{0.5pt}
  \drawline(120,80)(200,0)
  \drawline(120,80)(180,60)(200,0)

  \linethickness{0.3pt}

  \drawline(122,78)(126,78)
  \drawline(124,76)(132,76)
  \drawline(126,74)(138,74)
  \drawline(128,72)(144,72)
  \drawline(130,70)(150,70)
  \drawline(132,68)(156,68)
  \drawline(134,66)(162,66)
  \drawline(136,64)(168,64)
  \drawline(138,62)(174,62)
  \drawline(140,60)(180,60)
  \drawline(142,58)(180.67,58)
  \drawline(144,56)(181.33,56)
  \drawline(146,54)(182,54)
  \drawline(148,52)(182.67,52)
  \drawline(150,50)(183.33,50)
  \drawline(152,48)(184,48)
  \drawline(154,46)(184.67,46)
  \drawline(156,44)(185.33,44)
  \drawline(158,42)(186,42)
  \drawline(160,40)(186.67,40)
  \drawline(162,38)(187.33,38)
  \drawline(164,36)(188,36)
  \drawline(166,34)(188.67,34)
  \drawline(168,32)(189.33,32)
  \drawline(170,30)(190,30)
  \drawline(172,28)(190.67,28)
  \drawline(174,26)(191.33,26)
  \drawline(176,24)(192,24)
  \drawline(178,22)(192.67,22)
  \drawline(180,20)(193.33,20)
  \drawline(182,18)(194,18)
  \drawline(184,16)(194.67,16)
  \drawline(186,14)(195.33,14)
  \drawline(188,12)(196,12)
  \drawline(190,10)(196.67,10)
  \drawline(192,8)(197.33,8)
  \drawline(194,6)(198,6)
  \drawline(196,4)(198.67,4)
  \drawline(198,2)(199.33,2)

  \drawline(122,78)(122,79.33)
  \drawline(124,76)(124,78.67)
  \drawline(126,74)(126,78)
  \drawline(128,72)(128,77.33)
  \drawline(130,70)(130,76.67)
  \drawline(132,68)(132,76)
  \drawline(134,66)(134,75.33)
  \drawline(136,64)(136,74.67)
  \drawline(138,62)(138,74)
  \drawline(140,60)(140,73.33)
  \drawline(142,58)(142,72.67)
  \drawline(144,56)(144,72)
  \drawline(146,54)(146,71.33)
  \drawline(148,52)(148,70.67)
  \drawline(150,50)(150,70)
  \drawline(152,48)(152,69.33)
  \drawline(154,46)(154,68.67)
  \drawline(156,44)(156,68)
  \drawline(158,42)(158,67.33)
  \drawline(160,40)(160,66.67)
  \drawline(162,38)(162,66)
  \drawline(164,36)(164,65.33)
  \drawline(166,34)(166,64.67)
  \drawline(168,32)(168,64)
  \drawline(170,30)(170,63.33)
  \drawline(172,28)(172,62.67)
  \drawline(174,26)(174,62)
  \drawline(176,24)(176,61.33)
  \drawline(178,22)(178,60.67)
  \drawline(180,20)(180,60)
  \drawline(182,18)(182,54)
  \drawline(184,16)(184,48)
  \drawline(186,14)(186,42)
  \drawline(188,12)(188,36)
  \drawline(190,10)(190,30)
  \drawline(192,8)(192,24)
  \drawline(194,6)(194,18)
  \drawline(196,4)(196,12)
  \drawline(198,2)(198,6)

\end{picture}
\caption{Other operators analogous to $T$}\label{fig:Tgen}
\end{figure}

It is natural to consider more general operators than $T$, defined also as limits of random matrices or,
equivalently, in the approach we will take in this paper, 
by cutting a circular operator $Z$ using projections as in~\cite[\S4]{DH2}.
Some of these are pictured in Figure~\ref{fig:Tgen}, where again the shaded regions indicate weight $1$ and the unshaded
regions have weight $0$.
It is natural to ask whether these operators have nontrivial hyperinvariant subspaces.
The approach used in~\cite{DH3} for $T$ is not presently tenable, however;
while individual
$*$--moments for these operators can be calculated rather easily, a good general formula
is lacking;
moreover, such special relations between moments of $TT^*$ and $T^k(T^*)^k$ as mentioned above
are unlikely to be found in more general settings.

In this paper, we use another technique to exhibit nontrivial hyperinvariant subspaces for all operators
in a large class generalizing $T$, (including those pictured in Figure~\ref{fig:Tgen}).
This technique is an adaptation of one recently found by 
Foia\c s, Jung, Ko and Pearcy~\cite{FJKP}, which they applied to
certain quasinilpotent operators $Q$ in $\Bc(\HEu)$.
They consider spectral resolutions of $Q^k(Q^*)^k$ acting on vectors $x_0\in\HEu$.
Our modification, is, firstly, to take $Q$ in a II$_1$--factor $\Mcal$ and for $x_0$ to take the
trace vector in the standard representation of $\Mcal$,
and, secondly, to consider simultaneously a unital subalgebra $\Nc\subseteq\Mcal$ and the
conditional expectations of $Q^k(Q^*)^k$ onto $\Nc$ for positive integers $k$.

\begin{figure}[t]
\begin{picture}(200,200)(0,-30)

  \linethickness{0.7pt}
  \drawline(0,0)(0,160)(160,160)(160,0)(0,0) 

  \linethickness{0.5pt}
  \drawline(0,160)(160,0) 
  \put(50,50){$0$}
  \put(120,100){$*$}

  \drawline(25,135)(75,135)(75,85) 
  \put(56,116){$r$}

  \drawline(108,52)(160,52) 
  \put(130,35){$\le R$}

  \drawline(160,160)(190,160) 
  \drawline(187,162)(190,160)(187,158) 
  \put(175,165){$t$}
  \drawline(160,160)(160,165) 
  \put(154,168){$1$}
  \drawline(108,160)(108,165) 
  \put(102,168){$a$}
  \drawline(75,160)(75,165) 
  \put(69,168){$d$}
  \drawline(25,160)(25,165) 
  \put(19,168){$c$}

  \drawline(0,0)(0,-30) 
  \drawline(-2,-27)(0,-30)(2,-27) 
  \put(-15,-20){$s$}
  \drawline(0,0)(-5,0) 
  \put(-13,3){$1$}
  \drawline(0,52)(-5,52) 
  \put(-13,55){$a$}
  \drawline(0,85)(-5,85) 
  \put(-13,88){$d$}
  \drawline(0,135)(-5,135) 
  \put(-13,138){$c$}

\end{picture}
\caption{Conditions on $\eta$.}\label{fig:etacond}
\end{figure}
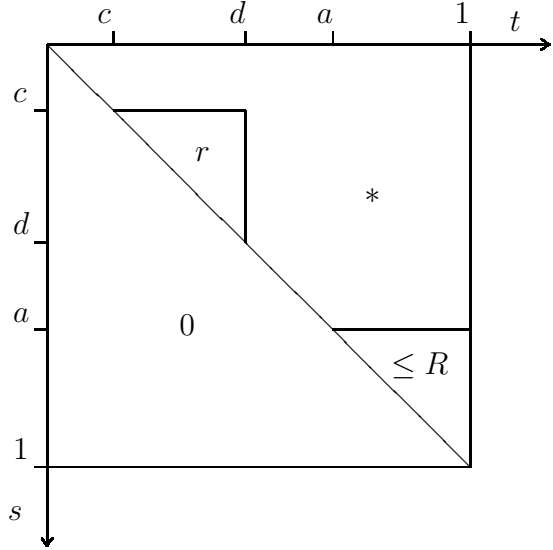

The class of operators we consider are certain $B$--circular
operators.
We introduce $B$--circular operators, which are
a special case of Speicher's $B$--Gaussian operators~\cite{Sp}.
Examples include the usual circular operator, Shlyakhtenko's generalized circular operators~\cite{Sh97},
the quasinilpotent DT--operator $T$ and the operators pictured in Figure~\ref{fig:Tgen}.
After proving some facts about $B$--circular operators,
we specialize to $B$--circular operators in tracial von Neumann algebras
when $B=L^\infty([0,1])$.
It turns out that these are the operators $z_\eta$, where $\eta$ is any finite
Borel measure on $[0,1]^2$ whose push--forwards ${\pi_i}_*\eta$ under the coordinate projections
$\pi_1,\pi_2:[0,1]^2\to[0,1]$ are absolutely continuous with respect to Lebesgue measure and
have essentially bounded Radon--Nikodym derivatives with respect to Lebesgue measure.
When $\eta$ is Lebesgue measure on $[0,1]^2$, then $z_\eta$ is the usual circular operator.
When $\eta$ is the restriction of Lebesgue measure to the upper triangle pictured in Figure~\ref{fig:T},
then $z_\eta$ is the quasinilpotent DT--operator $T$, while when $\eta$ is, for
example, the restriction
of Lebesgue measure to one of the shaded regions depicted in Figure~\ref{fig:Tgen}, then
$z_\eta$ is the corresponding generalization of $T$ described above.
We show that $z_\eta$ has a nontrivial hyperinvariant
subspace whenever the following three criteria hold:
\renewcommand{\labelenumi}{(\roman{enumi})}
\begin{enumerate}
\item $\eta$ is supported in the upper triangle $\{(s,t)\mid0\le s\le t\le1\}$;
\item for some $0<c<d$, the restriction of $\eta$ to $\{(s,t)\mid c\le s\le t\le d\}$
      is $r$ times Lebesgue measure, for some $r>0$;
\item for some $0<a<1$, the restriction of $\eta$ to $\{(s,t)\mid a\le s\le t\le1\}$
      is less than or equal to $R$ times Lebesgue measure, for some $R<\infty$.
\end{enumerate}
These conditions on $\eta$ are illustrated in Figure~\ref{fig:etacond}.
(Actually, some weaker conditions on $\eta$ suffice --- see Theorem~\ref{thm:zhis}
and Figure~\ref{fig:etagencond}.)

We now describe the contents of the rest of the paper.
In~\S\ref{sec:Wstar}, we show the well known fact that the projection
onto an $A$--hyperinvariant subspace belongs to the von Neumann algebra $\vN(A)$
generated by $A$.
We then show that, given an element $A$ of a W$^*$--algebra $\Mcal$,
the set of projections in $\Mcal$ that correspond to $A$--hyperinvariant subspaces
is independent of the normal $*$--representation of $\Mcal$.
The proof is technically straightforward, but the result is, we believe, conceptually valuable.
We also give some related examples.
In~\S\ref{sec:constr}, we prove a version of the construction of hyperinvariant subspaces
from~\cite{FJKP} applicable to certain operators in a tracial von Neumann algebra.
In~\S\ref{sec:Bcirc}, we introduce $B$--circular operators and prove several results
about them.
In~\S\ref{sec:hisp}, we use the method from~\S\ref{sec:constr} to construct
nontrivial hyperinvariant subspaces for the operators $z_\eta$ with $\eta$ satisfying
conditions (i)--(iii) above.
In~\S\ref{sec:inLF3}, we construct $z_\eta$ in $L(\Fb_3)$
when $\eta$ is absolutely continuous
with respect to Lebesgue measure on $[0,1]^2$,
using a method 
analogous to that of~\cite[\S4]{DH2}.
Finally, in~\S\ref{sec:quasinil}, we show that $z_\eta$ is quasinilpotent if $\eta$ is supported
on the upper triangle and is Lebesgue absolutely continuous with bounded Radon--Nikodym derivative
near the diagonal.

\medskip
\noindent
{\em Acknowledgements.}
The author is grateful to Carl Pearcy for providing him with an early copy of~\cite{FJKP},
to Ron Douglas and Carl Pearcy
for helpful discussions about hyperinvariant subspaces and to Lars Aagaard and Uffe Haagerup for discussions
pointing to the main idea of the proof of Proposition~\ref{prop:Rtransf}.

\section{Hyperinvariant subspaces of operators in W$^*$--algebras}
\label{sec:Wstar}

If $\HEu_0$ is a subspace of $\HEu$ and if $p:\HEu\to\HEu_0$ is the projection onto $\HEu_0$,
then $\HEu_0$ is $A$--invariant if and only if $Ap=pAp$.
(Throughout this paper, all projections will be assumed to be self--adjoint.)
We will say that a projection $p\in\Bc(\HEu)$ is $A$--invariant if $p\HEu$ is an $A$--invariant subspace.

Let $\Mcal\subseteq\Bc(\HEu)$ be a von Neumann algebra.
A subspace $\HEu_0\subseteq\HEu$ is said to be {\em affiliated} to $\Mcal$ if the projection $p:\HEu\to\HEu_0$
onto $\HEu_0$ belongs to $\Mcal$.
Of particular interest for an operator $A\in\Bc(\HEu)$ are $A$--invariant subspaces that are affiliated to
the von Neumann algebra $\vN(A)$ generated by $A$.

The following result is well known and easy to show.
\begin{prop}\label{prop:invN}
Given a Hilbert space $\HEu$ and an operator $A\in\Bc(\HEu)$, if $\HEu_0\subseteq\HEu$ is
an $A$--hyperinvariant subspace, then $\HEu_0$ is affiliated 
to the von Neumann algebra $\vN(A)$ generated by $A$ in $\Bc(\HEu)$.
\end{prop}
\begin{proof}
Let $p$ be the projection onto an $A$--hyperinvariant subspace.
Suppose $S$ is in the commutant of $\vN(A)$.
Then $S$ commutes with $A$, so $Sp=pSp$.
But also $S^*$ commutes with $A$, so $S^*p=pS^*p$ and $pSp=pS$.
Thus $pS=Sp$.
By von Neumann's double commutant theorem, $p\in\vN(A)$.
\end{proof}

However, there may be $A$--invariant subspaces that are affiliated with $\vN(A)$ but are
not $A$--hyperinvariant, as the following example shows (see also Examples~\ref{ex:6by6} and~\ref{ex:10by10}).
Indeed, this is not surprising,
because $\vN(A)$ incorporates information about how $A$ related to its adjoint $A^*$,
while the (abstract) lattice of hyperinvariant subspaces of $A$ is a similarity invariant.
\begin{example}\label{ex:A}
Let $A$ be $\left(\begin{smallmatrix}0&0&0\\0&0&1\\0&0&1\end{smallmatrix}\right)$
acting on a three--dimensional Hilbert space.
Then 
\[
\vN(A)=\left\{\left(\begin{smallmatrix}a&0&0\\0&b_{11}&b_{12}\\0&b_{21}&b_{22}\end{smallmatrix}\right)\bigg|a,b_{ij}\in\Cpx\right\}
\]
and $p=\left(\begin{smallmatrix}0&0&0\\0&1&0\\0&0&1\end{smallmatrix}\right)$
is the projection onto a subspace that is $A$--invariant and affiliated to $\vN(A)$, but not $B$--invariant,
where $B=\left(\begin{smallmatrix}0&1&-1\\0&0&0\\0&0&0\end{smallmatrix}\right)$.
Note that $B$ commutes with $A$, and thus the range of $p$ is not $A$--hyperinvariant.
\end{example}

It is in any case natural to ask,
when $p\in\vN(A)$ is a projection onto an $A$--hyperinvariant subspace
and when $\pi:\vN(A)\to\Bc(\KEu)$ is a normal, faithful $*$--homomorphism,
whether the range of $\pi(p)$ in $\KEu$ must be a $\pi(A)$--hyperinvariant subspace.
In other words, 
given an abstract W$^*$--algebra $\Mcal$ and an element $A\in\Mcal$,
are the projections onto $A$--hyperinvariant subspaces the same for all representations of $\Mcal$?

As is seen below in Theorem~\ref{thm:pip},
an affirmative answer to the above question follows readily from the classical result
that every normal, faithful $*$--homomorphism
of a von Neumann algebra is an amplification followed by an induction.

\begin{lemma}\label{lem:amp}
Let $\HEu$ and $\KEu$ be nonzero Hilbert spaces and let $\HEu_0$ be a subspace of $\HEu$.
Take $A\in\Bc(\HEu)$ and consider the operator $A\otimes I\in\Bc(\HEu\otimes\KEu)$.
Then $\HEu_0$ is $A$--hyperinvariant if and only if $\HEu_0\otimes\KEu$ is $(A\otimes I)$--hyperinvariant.
\end{lemma}
\begin{proof}
Suppose $\HEu_0\otimes\KEu$ is $(A\otimes I)$--hyperinvariant.
If $\xi\in\HEu_0$, $S\in\Bc(\HEu)$ and $AS=SA$, then $(S\otimes I)$ commutes with $(A\otimes I)$.
Fixing any $\eta\in \KEu$, we have $(S\xi)\otimes\eta=(S\otimes I)(\xi\otimes\eta)\in\HEu_0\otimes\KEu$
and, therefore, $S\xi\in\HEu_0$.
Thus, $\HEu_0$ is $A$--hyperinvariant.

On the other hand, suppose $\HEu_0$ is $A$--hyperinvariant.
If $\eta\in\KEu$, let $V_\eta:\HEu\to\HEu\otimes\KEu$ be the map $V_\eta(\xi)=\xi\otimes\eta$.
Then for every $B\in\Bc(\HEu)$, we have $V_\eta B=(B\otimes I)V_\eta$.
Suppose $X\in\Bc(\HEu\otimes\KEu)$ commutes with $A\otimes I$.
Given $\eta_1,\eta_2\in\KEu$, we have
\[
V_{\eta_2}^*XV_{\eta_1}A=V_{\eta_2}^*X(A\otimes I)V_{\eta_1}
=AV_{\eta_2}^*XV_{\eta_1},
\]
and we deduce $V_{\eta_2}^*XV_{\eta_1}\HEu_0\subseteq\HEu_0$.
Consequently, $X(\HEu_0\otimes\KEu)\subseteq\HEu_0\otimes\KEu$, and $\HEu_0\otimes\KEu$ is $(A\otimes I)$--hyperinvariant.
\end{proof}

\begin{lemma}\label{lem:ind}
Suppose $A\in\Bc(\HEu)$ and $\HEu_0\subseteq\HEu$ is an $A$--hyperinvariant subspace.
Let $P_0:\HEu\to\HEu_0$ be the projection onto $\HEu_0$ and suppose $E\in\Bc(\HEu)$
is a projection that commutes with $A$ and with $P_0$.
Let $A_E$ denote the operator in $\Bc(E\HEu)$ obtained by restricting $A$ to $E\HEu$.
Then $E\HEu\cap\HEu_0$ is $A_E$--hyperinvariant.
\end{lemma}
\begin{proof}
Let $S\in\Bc(E\HEu)$ commute with $A_E$.
Let $T=SE\in\Bc(\HEu)$.
Then $AT=TA$, so $T\HEu_0\subseteq\HEu_0$.
But $T\HEu\subseteq E\HEu$, so $S(E\HEu\cap\HEu_0)=T(\HEu_0)\subseteq E\HEu\cap\HEu_0$.
Therefore, $E\HEu\cap\HEu_0$ is $A_E$--hyperinvariant.
\end{proof}

We let $\Proj(\Mcal)$ denote the set of all projections, i.e.\ self--adjoint idempotents, in a von Neumann
algebra $\Mcal$.
As promised, the following theorem allows us to speak of hyperinvaraint projections
of an element of a von Neumann algebra, independent of representation on Hilbert space.

\begin{thm}\label{thm:pip}
Let $\Mcal\subseteq\Bc(\HEu)$ be a von Neumann algebra, let $A\in\Mcal$, let $p\in\Proj(\Mcal)$
and suppose $p\HEu$ is a $A$--hyperinvariant.
Let $\pi:\Mcal\to\Bc(\HEu_\pi)$ be any normal, faithful $*$--representation of $\Mcal$.
Then $\pi(p)(\HEu_\pi)$ is $\pi(A)$--hyperinvariant.
\end{thm}
\begin{proof}
By~\cite[Ch.\ I, \S4, Thm.\ 3]{Dix}, there is a Hilbert space $\KEu$, a projection $E$ in the commutant of
$\Mcal\otimes I_\KEu$ in $\Bc(\HEu\otimes\KEu)$ and a unitary $U:\HEu_\pi\to E(\HEu\otimes\KEu)$ such that
$\pi(x)=U^*(x\otimes I_\KEu)EU$.
But then $\pi(p)\HEu_\pi=U^*E((p\HEu)\otimes\KEu)$.
By Lemma~\ref{lem:amp}, $p\HEu\otimes\KEu$ is $(A\otimes I_\KEu)$--hyperinvariant.
From Lemma~\ref{lem:ind}, we then obtain that $E(p\HEu\otimes\KEu)=E(\HEu\otimes\KEu)\cap(p\HEu\otimes\KEu)$ is
$E(A\otimes I_\KEu)$--hyperinvariant.
Therefore, $U^*(E(p\HEu\otimes\KEu))=\pi(p)\HEu_\pi$ is $U^*(E(A\otimes I_\KEu))U$--hyperinvariant,
i.e.\ is $\pi(A)$--hyperinvariant.
\end{proof}

\begin{defi}
Let $\Mcal$ be a W$^*$--algebra, let $A\in\Mcal$ and let $p\in\Proj(\Mcal)$.
We call $p$ an $A$--{\em hyperinvariant projection} if $\pi(p)\HEu_\pi$ is
a $\pi(A)$--hyperinvariant subspace for one (and then for all) normal, faithful
$*$--homomorphisms $\pi:\Mcal\to\Bc(\HEu_\pi)$.
\end{defi}

\begin{remark}
By a result~\cite[Cor.\ 1.5]{DP} of Douglas and Pearcy, which
utilizes work of Hoover~\cite{H}, if $\Mcal$ is a von Neumann algebra that can be written
as a direct sum $\Mcal=\Mcal_1\oplus\Mcal_2$ with $\Mcal_1$
a (nonzero) finite type~I von Neumann algebra,
and if $A\in\Mcal$ is not a scalar multiple of the identity,
then $A$ has a nontrivial hyperinvariant projection.
\end{remark}

In light of Theorem~\ref{thm:pip},
it stands to reason that there should be representation--independent descriptions
(whatever that may mean) of the $A$--hyperinvariant projections in $\vN(A)$.
In that light, it seems natural to ask the following question.

\begin{ques}\label{qn}
Let $A$ be an operator in Hilbert space such that the von Neumann algebra $\vN(A)$ it generates
is a factor not isomorphic to $\Cpx$.
If $p$ is a projection in $\vN(A)$ and if $p$ is $S$--invariant for every
element $S$ of $\vN(A)$ that commutes with $A$, is $p$ necessarily an $A$--hyperinvariant projection?
\end{ques}
The answer is negative if we do not require $\vN(A)$ be be a factor;
indeed, the projection $p$ from Example~\ref{ex:A} belongs to the center of $\vN(A)$,
but fails to be $A$--hyperinvariant.
However, as far as the author knows, Question~\ref{qn} is open,
(though of course if $\vN(A)$ is a factor of type~I, then  the answer is positive).

In any case, Examples~\ref{ex:6by6} and~\ref{ex:10by10} below show that even when $A$ generates a factor of type I or of type II$_1$,
there may be an $A$--invariant subspace affiliated to the factor that is not $A$--hyperinvariant.

We need a preparatory, elementary lemma about $n\times n$ matrices.
Let $\{e_{i,j}\mid1\le i,j\le n\}$ be a system of matrix units in $M_n(\Cpx)$.
\begin{lemma}\label{lem:irreduc}
Let $n,p\in\Nats$ with $p\ge2$, $n>2p$.
Let $b_1,\ldots,b_{n-p}$ be distinct, strictly positive numbers.
Then there is $\eps>0$ such that whenever $a_2,\ldots,a_p\in(0,\eps)$ and
\[
A=\sum_{k=1}^{n-p}b_ke_{k,k+p}+\sum_{k=2}^pa_ke_{1,p+k},
\]
the $*$--algebra generated by $A$ is all of $M_n(\Cpx)$.
\end{lemma}
\begin{proof}
We may write
\[
A=\left(\begin{smallmatrix}
0&0&\cdots&0&b_1&a_2&\cdots&a_p&0      &\cdots&0 \\
 & &      & &   &b_2 \\
 & &      & &   &   &\ddots \\
 & &      & &   &   &      &b_p \\
 & &      & &   &   &      &   &b_{p+1} \\
 & &      & &   &   &      &   &       &\ddots \\
 & &      & &   &   &      &   &       &      & b_{n-p}
\end{smallmatrix}\right),
\]
where the omitted entries are zero.
Let $\Afr$ denote the $*$--algebra generated by $A$.
Let
\[
B=\sum_{k=1}^pb_ke_{k,k+p}+\sum_{k=2}^pa_ke_{1,p+k},
\]
so that
\[
A=B+\sum_{k=p+1}^{n-p}b_ke_{k,k+p}.
\]
Then
\[
AA^*=BB^*+\sum_{k=p+1}^{n-p}b_k^2e_{k,k},
\]
while
\[
BB^*=\sum_{k=1}^pb_k^2e_{k,k}+\big(\sum_{k=2}^pa_k^2\big)e_{1,1}
+\sum_{k=2}^pb_ka_k(e_{1,k}+e_{k,1}).
\]
By choosing $\eps$ sufficiently small, the nonzero eigenvalues of $BB^*$ can be forced
to be arbitrarily close to $b_1^2,b_2^2,\ldots,b_p^2$, respectively.
Then we obtain
\begin{equation}\label{eq:Afr1}
e_{k,k}\in\Afr,\qquad(p+1\le k\le n-p)
\end{equation}
by taking spectral projections of $AA^*$.
We have $Ae_{p+1,p+1}=b_1e_{1,p+1}\in\Afr$, so 
\begin{equation}\label{eq:Afr2}
e_{1,1},e_{1,p+1}\in\Afr.
\end{equation}

From 
\begin{equation*}
(A-e_{1,1}A)e_{k+p,k+p}=b_ke_{k,k+p},\qquad(2\le k\le n-p),
\end{equation*}
together with~\eqref{eq:Afr1} and~\eqref{eq:Afr2},
we get
\begin{equation}\label{eq:Afr4}
e_{k,k+p}\in\Afr,\qquad(1\le k\le n-2p).
\end{equation}
Combined with~\eqref{eq:Afr1} this yields
\begin{equation*}
e_{k,k}\in\Afr,\qquad(1\le k\le n-p).
\end{equation*}
Using
\[
e_{k,k}A=b_ke_{k,k+p},\qquad(2\le k\le n-p),
\]
combined with~\eqref{eq:Afr4}, we now get
\begin{equation}\label{eq:Afr6}
e_{k,k+p}\in\Afr,\qquad(1\le k\le n-p).
\end{equation}
From
\begin{equation*}
Ae_{k+p,k+p}-b_ke_{k,k+p}=a_ke_{1,k+p},\qquad(2\le k\le p)
\end{equation*}
together with~\eqref{eq:Afr2}, we get
\begin{equation*}
e_{1,k+p}\in\Afr,\qquad(1\le k\le p)
\end{equation*}
and thus also
\[
e_{1,k}=e_{1,k+p}e_{k+p,k}\in\Afr,\qquad(2\le k\le p)
\]
and, using~\eqref{eq:Afr6},
\[
e_{1,k+2p}=e_{1,k+p}e_{k+p,k+2p}\in\Afr,\qquad(1\le k\le p).
\]
Continuing as long as possible, we get
\begin{align*}
e_{1,k+3p}&=e_{1,k+2p}e_{k+2p,k+3p}\in\Afr,\qquad(1\le k\le\min(p,n-3p)), \\
e_{1,k+4p}&=e_{1,k+3p}e_{k+3p,k+4p}\in\Afr,\qquad(1\le k\le\min(p,n-4p)), \\
\quad\vdots&
\end{align*}
yielding
\[
e_{1,j}\in\Afr,\qquad(1\le j\le n)
\]
and $\Afr=M_n(\Cpx)$.
\end{proof}

\begin{example}\label{ex:6by6}
We will find an operator $A$, generating a finite type~I factor
and having an invaraint subspace affiliated to the factor that is not,
however, $A$--hyperinvariant.

Let $a>0$ and consider the upper--triangular $6\times 6$ matrix
\[
A=\left(
\begin{smallmatrix}
0 & 0 & 1 & a &   &   \\
  & 0 & 0 & 2 &   &   \\
  &   & 0 & 0 & 3 &   \\
  &   &   & 0 & 0 & 4 \\
  &   &   &   & 0 & 0 \\
  &   &   &   &   & 0
\end{smallmatrix}\right),
\]
where the omitted entries are all zero.
By Lemma~\ref{lem:irreduc},
for $a$ sufficiently small we have $\vN(A)=M_6(\Cpx)$.
But taking the Jordan canonical form, $A$ is similar to
\[
B=\left(
\begin{smallmatrix}
0 & 1 &   &   &   &   \\
  & 0 & 1 &   &   &   \\
  &   & 0 & 0 &   &   \\
  &   &   & 0 & 1 &   \\
  &   &   &   & 0 & 1 \\
  &   &   &   &   & 0
\end{smallmatrix}\right),
\]
say $A=SBS^{-1}$.
Thus, the subspace that is the range of the idempotent operator $S\diag(1,1,1,0,0,0)S^{-1}$
is $A$--invariant and is affiliated
to $\vN(A)$, but is not $A$--hyperinvariant, because it is not invariant under
$SCS^{-1}$, where $C=\left(\begin{smallmatrix}0_3 & I_3 \\ I_3 & 0_3\end{smallmatrix}\right)$.
\end{example}

\begin{example}\label{ex:10by10}
We will find an operator $A$, generating a type~II$_1$ factor
and having an invaraint subspace affiliated to the factor that is not,
however, $A$--hyperinvariant.

Let $a>0$ and consider the upper--triangular $10\times10$ matrix
\[
F=\left(
\begin{smallmatrix}
0 & 0 & 1 & a & a & a \\
  & 0 & 0 & 2         \\
  &   & 0 & 0 & 3     \\
  &   &   & 0 & 0 & 4 \\
  &   &   &   & 0 & 0 & 5 \\
  &   &   &   &   & 0 & 0 & 6 \\
  &   &   &   &   &   & 0 & 0 & 7 \\
  &   &   &   &   &   &   & 0 & 0 & 8 \\
  &   &   &   &   &   &   &   & 0 & 0 \\
  &   &   &   &   &   &   &   &   & 0
\end{smallmatrix}\right),
\]
where again the omitted entries are all zero.
Then
\[
F^2=\left(
\begin{smallmatrix}
0 & 0 & 0 & 0 & 3 & 4a & 5a & 6a   \\
\phantom{48}  & 0 & 0 & 0 & 0 & 8          \\
  & \phantom{48}  & 0 & 0 & 0 & 0 & 15     \\
  &   & \phantom{48}  & 0 & 0 & 0 & 0  & 24  \\
  &   &   & \phantom{48}  & 0 & 0 & 0  & 0  & 35  \\
  &   &   &   & \phantom{48} & 0 & 0  & 0  & 0  & 48 \\
  &   &   &   &   &   & 0  & 0  & 0  & 0  \\
  &   &   &   &   &   &    & 0  & 0  & 0  \\
  &   &   &   &   &   &    &    & 0  & 0  \\
  &   &   &   &   &   &    &    &    & 0
\end{smallmatrix}\right).
\]
By Lemma~\ref{lem:irreduc}, for sufficienatly small $a$ we have
\begin{equation}\label{eq:F2}
\vN(F^2)=M_{10}(\Cpx).
\end{equation}

Suppose a II$_1$--factor $\Mcal$ is generated by $\{u,b\}$, where $u$ is a unitary satisfying $u^2=1$
and where $b\ge0$ and $b$ has spectrum in $[1,1+\eps]$ for some $\eps>0$ to be determined later.
For example, (see~\cite{D}), the interpolated free group factors $L(\Fb_t)$ for any $t\in(1,\frac32]$
have generators with these properties.
Let $x=ub$ and consider $A=F\otimes x\in M_{10}(\Cpx)\otimes\Mcal$.
Let $F^*F=\sum_{i=1}^n\lambda_iP_i$, where $P_1,\ldots,P_n$ are orthogonal projections and $\lambda_1,\ldots,\lambda_n$
are the distinct, nonzero eigenvalues of $F^*F$.
Then
\[
A^*A=F^*F\otimes b^2=\sum_{i=1}^n\lambda_iP_i\otimes b^2.
\]
If $\eps$ is small enough, then by taking spectral projections we get
$Q\otimes b^2\in\vN(A)$, where $Q=\sum_{i=1}^nP_i$.
Therefore,
$Q\otimes b^{-1}\in\vN(A)$ and
\[
F\otimes u=A(Q\otimes b^{-1})\in\vN(A).
\]
But $(F\otimes u)^2=F^2\otimes1$.
From~\eqref{eq:F2}, we get $M_{10}(\Cpx)\otimes1\subseteq\vN(A)$.
Consequently, $1\otimes ub\in\vN(A)$, and $A$ generates the II$_1$--factor $M_{10}(\Cpx)\otimes\Mcal$.

However, $F$ is similar in $M_{10}(\Cpx)$ to its Jordan canonical form
\[
G=\left(
\begin{smallmatrix}
0 & 1  \\
  & 0 & 1  \\
  &   & 0 & 1  \\
  &  &   & 0 & 1  \\
  &  &   &   & 0 & 0  \\
  &  &   &   &  & 0 & 1  \\
  &  &   &   &  &   & 0 & 1 \\
  &  &   &   &   &  &   & 0 & 1 \\
  &  &   &   &   &   &  &   & 0 & 1 \\
  &  &   &   &  &   &   &   &   & 0
\end{smallmatrix}\right),
\]
so $A=F\otimes x$ is similar in $\vN(A)=M_{10}(\Cpx)\otimes\Mcal$ to $G\otimes x$.
Arguing as in Example~\ref{ex:6by6}, we find a subspace that is $A$--invariant but not $A$--hyperinvariant and whose projection
lies in $M_{10}(\Cpx)\otimes1\subseteq\vN(A)$.
\end{example}

\section{A construction of hyperinvariant subspaces}
\label{sec:constr}

Foia\c s, Jung, Ko and Pearcy~\cite{FJKP} recently found a technique
that constructs
nontrivial hyperinvariant subspaces of some operators on Hilbert space.
In this section, we adapt their method so that it will apply to certain operators in tracial von Neumann algebras.

Let $\Mcal$ be a W$^*$--algebra having a normal, faithful, tracial state $\tau$.
We will consider the singular numbers of operators $a\in\Mcal$ with respect to $\tau$, which were
treated by Fack in~\cite{F} and by Fack and Kosaki in~\cite{FK}.
Thus, for $t\in[0,1]$, the $t$-th singular number of $a$ is
\begin{equation}\label{eq:sta}
s_t(a)=\inf\{\|a(1-p)\|\mid p\in\Proj(\Mcal),\,\tau(p)\le t\}.
\end{equation}
Of course, the singular numbers are highly dependent on the choice of trace $\tau$.
We may write $s_t(a;\tau)$ instead of $s_t(a)$, in order to avoid any confusion.
By~\cite[2.2]{FK}, we have
\begin{equation}\label{eq:FK}
s_t(a)=\inf\{\lambda\ge0\mid\tau(1_{(\lambda,\infty)}(|a|))\le t\},
\end{equation}
and the infimum is attained.
Here, $1_{(\lambda,\infty)}(|a|)$ denotes the Borel functional calculus,
so for $B\subseteq[0,\infty)$ and $x\in\Mcal$, $x\ge0$, $1_B(x)$ denotes
the spectral projection for $x$ corresponding to the set $B$.

Let $\Mcal$ be represented on the Hilbert space $L^2(\Mcal,\tau)$ via the Gelfand--Naimark--Segal construction.
Given $x\in\Mcal$, we will let $\xh$ denote the corresponding element of $L^2(\Mcal,\tau)$.
Suppose $\Nc\subseteq\Mcal$ is a unital W$^*$--subalgebra and let $\Ec:\Mcal\to\Nc$ be the $\tau$--preserving
conditional expectation onto $\Nc$.
As is well known, $\Ec$ is obtained by compressing with respect to the projection 
$e:L^2(\Mcal,\tau)\to L^2(\Nc,\tau\restrict_\Nc)$
onto the subspace of $L^2(\Mcal,\tau)$ that is identified with $L^2(\Nc,\tau\restrict_\Nc)$ by the inclusion $\Nc\subseteq\Mcal$.
In particular 
\begin{equation}\label{eq:e}
exe=e\Ec(x)=\Ec(x)e\qquad(x\in\Mcal).
\end{equation}

\begin{thm}\label{thm:FJKP}
Let $b\in\Mcal$.
Suppose there are integers $p\ge0$ and $1\le n(1)<n(2)<\cdots$ and there are $\theta\in(0,1)$ and $\mu_k\ge0$, ($k\in\Nats$),
so that
\begin{equation}\label{eq:ratio}
\lim_{k\to\infty}\frac{\mu_k}{s_\theta(b^{n(k)})^2}=0
\end{equation}
and there are vectors $\zeta_k\in L^2(\Nc,\tau\restrict_\Nc)$ such that 
\begin{equation}\label{eq:etak}
\zeta_k=1_{[0,\mu_k]}(\Ec(b^{n(k)+p}(b^*)^{n(k)+p}))\zeta_k
\end{equation}
and $\zeta_k$ converges with respect to the Hilbert--space norm on $L^2(\Nc,\tau\restrict_\Nc)$ to a nonzero
vector $\zeta\in L^2(\Nc,\tau\restrict_\Nc)$ as $k\to\infty$.
Then $b$ has a nontrivial, hyperinvariant subspace.
\end{thm}
\begin{proof}
For every $k\in\Nats$, take a sequence $(c_{k,j})_{j=1}^\infty$ in $\Nc$ so that $\widehat{c_{k,j}}$
converges to $\zeta_k$ as $j\to\infty$.
We may without loss of generality replace $\Mcal$ by the smallest von Neumann algebra
such that $b\in\Mcal$, all $c_{k,j}\in\Mcal$ and $\Ec(\Mcal)\subseteq\Mcal$.
Then $\Mcal$ is countably generated and $L^2(\Mcal,\tau)$ is separable.
Let $\lambda_k=s_\theta(b^{n(k)})^2=s_\theta(b^{n(k)}(b^*)^{n(k)})$.
For $n\in\Nats$, let $E_n$ be the projection--valued spectral measure of $b^n(b^*)^n$;
let
\[
x_k=\int_{[\lambda_k,\infty)}\frac1tdE_{n(k)}(t)\in\Mcal
\]
and $y_k=(b^*)^{n(k)}x_k$.
Then $b^{n(k)}y_k=E_{n(k)}([\lambda_k,\infty))$ and
\[
\|b^{n(k)}\widehat{y_k}\|_2^2=\langle E_{n(k)}([\lambda_k,\infty))\widehat{\:},\oneh\rangle=\tau(E_{n(k)}([\lambda_k,\infty)).
\]
From~\eqref{eq:FK}, we have $\tau(E_{n(k)}((\lambda,\infty))>\theta$ whenever $\lambda<\lambda_k$, so
\begin{equation}\label{eq:theta1}
\theta\le\inf_{\lambda<\lambda_k}\tau(E_{n(k)}((\lambda,\infty)))=\tau(E_{n(k)}([\lambda_k,\infty)))\le1.
\end{equation}
Since $\|b^{n(k)}\widehat{y_k}\|_2$ stays bounded as $k\to\infty$, by passing to a subsequence, if necessary,
we may without loss of generality assume $b^{n(k)}\widehat{y_k}$ converges in the weak topology to a vector $\xi\in L^2(\Mcal,\tau)$
as $k\to\infty$.
From~\eqref{eq:theta1}, we then have $\langle\xi,\oneh\rangle\ge\theta$, so $\xi\ne0$.
Moreover, we have
\begin{equation}\label{eq:yk}
\begin{aligned}
\|\widehat{y_k}\|_2^2&=\tau(x_kb^{n(k)}(b^*)^{n(k)}x_k) \\
&=\tau\bigg(b^{n(k)}(b^*)^{n(k)}\int_{[\lambda_k,\infty)}\frac1{t^2}dE_{n(k)}(t)\bigg)
=\tau(x_k)\le\frac1{\lambda_k}.
\end{aligned}
\end{equation}

Since $\zeta_k\in L^2(\Nc,\tau\restrict_\Nc)$, using~\eqref{eq:e} and~\eqref{eq:etak} we have
\begin{equation}\label{eq:b*eta}
\begin{aligned}
\|(b^*)^{n(k)+p}\zeta_k\|_2^2
&=\langle eb^{n(k)+p}(b^*)^{n(k)+p}e\zeta_k,\zeta_k\rangle
=\langle\Ec(b^{n(k)+p}(b^*)^{n(k)+p})\zeta_k,\zeta_k\rangle \\
&=\langle\Ec(b^{n(k)+p}(b^*)^{n(k)+p})1_{[0,\mu_k]}(\Ec(b^{n(k)+p}(b^*)^{n(k)+p}))\zeta_k,\zeta_k\rangle \\
&\le\mu_k\|\zeta_k\|_2^2.
\end{aligned}
\end{equation}

If $S\in\Bc(L^2(\Mcal,\tau))$ and if $S$ commutes with $b$, then we have
\[
\langle S\xi,(b^*)^p\zeta\rangle=\lim_{k\to\infty}\langle Sb^{n(k)}\widehat{y_k},(b^*)^p\zeta_k\rangle
=\lim_{k\to\infty}\langle S\widehat{y_k},(b^*)^{n(k)+p}\zeta_k\rangle.
\]
But from~\eqref{eq:yk} and~\eqref{eq:b*eta},
\[
|\langle S\widehat{y_k},(b^*)^{n(k)+p}\zeta_k\rangle|
\le\|S\|\|\zeta_k\|_2\sqrt{\frac{\mu_k}{\lambda_k}}.
\]
By hypothesis, this upper bound tends to zero as $k\to\infty$.
Therefore, we have
\begin{equation}\label{eq:Sxi}
\langle S\xi,(b^*)^p\zeta\rangle=0.
\end{equation}
Clearly,
\[
\Vc:=\overline{\{S\xi\mid S\in\Bc(L^2(\Mcal,\tau)),\, Sb=bS\}}
\]
is a nonzero $b$--hyperinvariant subspace.
If $(b^*)^p\zeta\ne0$, then by~\eqref{eq:Sxi}, $\Vc$ is nontrivial.
If, on the other hand, $(b^*)^p\zeta=0$, then $b$ has a nonzero cokernel.
Since $b$ is not the zero operator, it follows that
$\overline{b(L^2(\Mcal,\tau))}$ is a nontrivial $b$--hyperinvariant subspace.
\end{proof}

We will make use of the following well known result in application of Theorem~\ref{thm:FJKP}.

\begin{lemma}\label{lem:s}
Let $\Mcal$ be a von Neumann with normal, faithful, tracial state $\tau$,
let $a\in\Mcal$ and $q\in\Proj(\Mcal)$.
If $0<\theta<\tau(q)$, then
\begin{equation}\label{eq:s}
s_\theta(a;\tau)\ge s_{\frac\theta{\tau(q)}}(qaq;\tau(q)^{-1}\tau\restrict_{q\Mcal q}).
\end{equation}
\end{lemma}
\begin{proof}
Suppose $p\in\Proj(\Mcal)$, $\tau(p)\le\theta$.
Then
\[
\tau(q\wedge(1-p))\ge\tau(q)+\tau(1-p)-1\ge\tau(q)-\theta,
\]
so
\[
\frac{\tau(q-q\wedge(1-p))}{\tau(q)}\le\frac\theta{\tau(q)},
\]
and
\[
\|a(1-p)\|\ge\|a(q\wedge(1-p))\|\ge\|qaq(q\wedge(1-p))\|.
\]
This implies~\eqref{eq:s} directly from the definition~\eqref{eq:sta}.
\end{proof}

\section{$B$--circular elements}
\label{sec:Bcirc}

\begin{defi}\label{def:Bprob}
Let $B$ be a unital $*$--algebra over $\Cpx$.
\renewcommand{\labelenumi}{(\roman{enumi})}
\begin{enumerate}
\item A {\em $B$--valued $*$--noncommutative probability space} is a pair $(A,E)$,
where $A$ is a unital $*$--algebra containing $B$ as a unital $*$--subalgebra
(which makes $A$ into a bimodule over $B$)
and where $E:A\to B$ is
a $B$--bimodule map satisfying $E(b)=b$ for all $b\in B$.
\item We say $(A,E)$ is a {\em $B$--valued Banach $*$--noncommutative probability space} if, in addition,
$A$ is a unital Banach $*$--algebra, $B$ is a closed subalgebra of $A$ and $E$ is bounded.
\item We say $(A,E)$ is a {\em $B$--valued C$^*$--noncommutative probability space} if, in addition,
$A$ is a unital C$^*$--algebra,
$B$ a C$^*$--subalgebra of $A$ and $E$ is a projection of norm $1$ onto $B$.
(It follows from~\cite{Tom} that then $E$ is positive and a $B$--bimodule map.)
\item We say $(A,E)$ is a {\em $B$--valued W$^*$--noncommutative probability space} if, in addition,
$A$ is a unital W$^*$--algebra,
$B$ a W$^*$--subalgebra of $A$ and $E$ is a normal projection of norm $1$ onto $B$.
\end{enumerate}
\end{defi}

\begin{defi}\label{def:Bcirc}
Let $(A,E)$ be a $B$--valued $*$--noncommutative probability space.
Let $\alpha:B\to B$ and $\beta:B\to B$ be $\Cpx$--linear maps.
A {\em $B$--circular element} with covariance $(\alpha,\beta)$ is an element $z\in A$
such that the distribution of the pair $(z,z^*)$ is $B$--Gaussian in the sense
of~\cite[Def.\ 4.2.3]{Sp}, with covariance determined by
\begin{gather*}
E(z^*bz)=\alpha(b) \\
E(zbz^*)=\beta(b) \\
E(zbz)=E(z^*bz^*)=E(z)=E(z^*)=0
\end{gather*}
for all $b\in B$.
In the case that $(A,E)$ is a $B$--valued C$^*$--noncommutative probability space, we may call $z$
a {\em $B$--circular operator}.
\end{defi}

\begin{examples}\label{ex:Bcirc}
\renewcommand{\labelenumi}{(\roman{enumi})}
\begin{enumerate}
\item A usual circular operator $z$ with $\tau(z^*z)=r$ is, in the notation of Definition~\ref{def:Bcirc},
a $\Cpx$--circular element with covariance $(r,r)$, where here $r$ denotes multiplication by $r$ on $\Cpx$.
\item The generalized circular elements $\ell_2+\sqrt\lambda\ell_1^*$, ($0\le\lambda\le 1$), considered in~\cite{Sh97},
are $\Cpx$--circular with covariance $(\lambda,1)$, where again the scalars indicate operations of multiplication on $\Cpx$.
\item A $\DT(\delta_0,c)$ operator, considered in~\cite{DH2}, is $L^\infty([0,1])$--circular, with covariance
$(\alpha,\beta)$, where
\begin{align*}
\alpha(f)(x)&=c^2\int_0^xf(t)dt \\
\beta(f)(x)&=c^2\int_x^1f(t)dt.
\end{align*}
This was shown in the appendix to~\cite{DH3}.
\end{enumerate}
\end{examples}

The $B$--valued $*$--moments of a $B$--circular operator can be calculated using Speicher's free cummulant
calculus~\cite{Sp}.
This is amounts to the nested evaluation described by \'Sniady in~\cite[\S4.2]{Sn}.
This technique is reviewed below, using the notation $\pi\{\cdots\}$ for the bracketing of a noncrossing pair partition $\pi$
with the multiplicative function of free cummulants as in~\cite{NSS}.

\begin{remark}\label{rmk:mom}
With $z$ a $B$--circular operator as above, let $n\in\Nats$, $s(1),\ldots,s(n)\in\{1,*\}$ and $b_1,\ldots,b_n\in B$.
Then
\begin{equation}\label{eq:zmom}
E(z^{s(1)}b_1z^{s(2)}b_2\cdots z^{s(n)}b_n)
=\sum_{\pi\in\NC_2(n)}\pi\{z^{s(1)}b_1,\ldots,z^{s(n)}b_n\}
\end{equation}
where the sum is over all non--crossing pair partitions $\pi$ of $\{1,\ldots,n\}$ and where
the quantity
\begin{equation}\label{eq:pibrak}
\pi\{z^{s(1)}b_1,\ldots,z^{s(n)}b_n\},
\end{equation}
which is the bracketing of the cummulants of the pair $(z,z^*)$,
is evaluated as described below.
In particular, the $*$--moment~\eqref{eq:zmom} vanishes if $n$ is odd; so let us assume $n$ is even.
Let 
\begin{equation}\label{eq:pi}
\pi=\{\{i_1,j_1\},\ldots,\{i_{n/2},j_{n/2}\}\}.
\end{equation}
Then the quantity~\eqref{eq:pibrak} vanishes unless $s(i_p)\ne s(j_p)$ for all $p\in\{1,\ldots,n/2\}$,
i.e.\ unless $\pi$ pairs only $z$ with $z^*$.
Therefore, the $*$--moment~\eqref{eq:zmom} vanishes if the number of $j$ such that $s(j)=*$ differs from
the number of $j$ such that $s(j)=1$.
The quantity~\eqref{eq:pibrak} is evaluated as follows.
Without loss of generality take $i_1=1$ in~\eqref{eq:pi}.
Then
\begin{equation}\label{eq:piinduc}
\pi\{z^{s(1)}b_1,\ldots,z^{s(n)}b_n\}
=\left\{
\begin{alignedat}{2}
&\alpha\big(b_1(\pit'\{z^{s(2)}b_2,\ldots,z^{s(j_1-1)}b_{j_1-1}\})\big) \\
&\quad b_{j_1}(\pit''\{z^{s(j_1+1)}b_{j_1+1},\ldots,z^{s(n)}b_n\})
&\quad\text{if }s(1)=*, \\
&\beta\big(b_1(\pit'\{z^{s(2)}b_2,\ldots,z^{s(j_1-1)}b_{j_1-1}\})\big) \\
&\quad b_{j_1}(\pit''\{z^{s(j_1+1)}b_{j_1+1},\ldots,z^{s(n)}b_n\})
&\text{if }s(1)=1,
\end{alignedat}
\right.
\end{equation}
where $\pit'$ is the restriction of $\pi$ to $\{2,\ldots,j_1-1\}$,
renumbered by left translation to become an element of $\NC_2(j_1-2)$,
while $\pit''$
is the restriction of $\pi$ to $\{j_1+1,\ldots,n\}$, renumbered by translation to become
an element of $\NC_2(n-j_1)$, and where we set $\pit'\{z^{s(2)}b_2,\ldots,z^{s(j_1-1)}b_{j_1-1}\}$
to be $1$ if $j_1=2$ and $\pit''\{z^{s(j_1+1)}b_{j_1+1},\ldots,z^{s(n)}b_n\}$ to be $1$ if $j_1=n$.
\end{remark}

For example, if $n=6$ and $\pi=\{\{1,4\},\{2,3\},\{5,6\}\}$, then
\begin{align*}
\pi\{zb_1,z^*b_2,zb_3,z^*b_4,z^*b_5,zb_6\}
&=\beta\big(b_1(\pit'\{z^*b_2,zb_3\})\big)b_4\alpha(b_5)b_6 \\
&=\beta(b_1\alpha(b_2)b_3)b_4\alpha(b_5)b_6,
\end{align*}
with $\pit'=\{\{1,2\}\}$.

The following basic properties are special instances of Speicher's results~\cite{Sp}.

\begin{prop}\label{prop:properties}
Let $(A,E)$ be a $B$--valued $*$--noncommutative probability space and let $z$ and $z'$ be $B$--circular
elements in $(A,E)$ with covariances $(\alpha,\beta)$ and $(\alpha',\beta')$, respectively.
Suppose $z$ and $z'$ are $*$--free over $B$ with respect to $E$.
Then:
\renewcommand{\labelenumi}{(\roman{enumi})}
\begin{enumerate}
\item $z^*$ is $B$--circular with covariance $(\beta,\alpha)$.
\item $z+z'$ is $B$--circular with covariance $(\alpha+\alpha',\beta+\beta')$.
\item Let $d\in B$; then $d^*zd$ is $B$--circular with covariance $(\alpha_d,\beta_d)$, where
\begin{align*}
\alpha_d(b)&=d^*\alpha(dbd^*)d \\
\beta_d(b)&=d^*\beta(dbd^*)d.
\end{align*}
\item Suppose $p$ is a self--adjoint idempotent in $B$;
then in the $pBp$--valued $*$--noncommutative probability space $(pAp,E\restrict_{pAp})$,
$pzp$ is $pBp$--circular with covariance $(\alphat_p,\betat_p)$, where
$\alphat_p,\betat_p:pBp\to pBp$ are given by
\begin{align*}
\alphat_p(b)&=p\alpha(b)p \\
\betat_p(b)&=p\beta(b)p.
\end{align*}
\end{enumerate}
\end{prop}
\begin{proof}
Part~(i) is clear from the cummulant calculus.
Part~(ii) follows from the additivity of free cummulants of $*$--free variables~\cite[Thm.\ 4.1.7]{Sp}.
Part~(iii) follows from~\cite[Prop.\ 4.1.10]{Sp}.
Part~(iv) follows from part~(iii).
\end{proof}

\begin{prop}\label{prop:Bsemicirc}
Let $B$ be a unital $*$--algebra, let $(A,E)$ be a $B$--valued $*$-probability space,
let $\alpha,\beta:B\to B$ be $\Cpx$--linear maps and let $z\in A$.
Then $z$ is a $B$--circular operator with covariance $(\alpha,\beta)$ if and only if 
\begin{equation}\label{eq:ReIm}
z=\frac{x_1+ix_2}{\sqrt 2},\qquad x_i^*=x_i,
\end{equation}
where in the notation of~\cite[Def.\ 4.2.3]{Sp},
the distribution of the pair $x_1,x_2$ is $B$--Gaussian with covariance determined by
\begin{align*}
E(x_1bx_1)&=(\alpha(b)+\beta(b))/2 \\
E(x_1bx_2)&=i(\beta(b)-\alpha(b))/2 \\
E(x_2bx_1)&=i(\alpha(b)-\beta(b))/2 \\
E(x_2bx_2)&=(\alpha(b)+\beta(b))/2 \\
E(x_1)&=E(x_2)=0.
\end{align*}
\end{prop}
\begin{proof}
From~\eqref{eq:ReIm} we have $x_1=\frac{z+z^*}{\sqrt2}$ and $x_2=\frac{z-z^*}{\sqrt2}$.
Now the remaining assersions follow from multilinearity of $B$--valued cummulants.
\end{proof}

Compare the following result to~\cite[Prop.\ 2.20]{Sh},
from which it follows in light of Proposition~\ref{prop:Bsemicirc}.

\begin{prop}\label{prop:tracial}
Let $(A,E)$ be a $B$--valued $*$--noncommutative probability space and let $z\in A$ be $B$--circular
with covariance $(\alpha,\beta)$.
Suppose $\tau:B\to B$ is a trace on $B$.
Then the restriction of $\tau\circ E$ to the $*$--algebra generated by $B\cup\{z\}$ is a trace
if and only if $\tau(\alpha(b)c)=\tau(\beta(c)b)$ for all $b,c\in B$.
\end{prop}

\begin{notation}
Let $(A,E)$ be a $B$--valued Banach $*$--noncommutative probability space.
Given $x\in A$, for $b\in B$ with $\|b\|$ sufficiently small,
we set
\[
\Gt_x(b)=\sum_{n=0}^\infty E(b(xb)^n).
\]
Note that $\Gt_x$ is related to the Cauchy transform $G_x(b)=E((b-x)^{-1})$, as it appears
for example in~\cite{V} or~\cite{VDN},
by
$G_x(b)=\Gt_x(b^{-1})$ for $b\in B$ invertible with $\|b^{-1}\|$ sufficiently small.
\end{notation}

At the heart of the proof of the following result
is a scheme for finding the generating function of the Catalan numbers, when
we recall that the Catalan number $C_n=\frac1{n+1}\binom{2n}n$
is the number of non--crossing pair partitions of $\{1,\ldots,2n\}$.
The author is indebted to Lars Aagaard and Uffe Haagerup for discussions
of this method for the quasinilpotent DT--operator.

\begin{prop}\label{prop:Rtransf}
Suppose $(A,E)$ is a $B$--valued Banach $*$--noncommutative probability space
and $z\in A$ is a $B$--circular operator with covariance $(\alpha,\beta)$.
Then for $b,c\in B$ with $\|b\|\|c\|$ sufficiently small, we have
\begin{align}
\Gt_{z^*cz}(b)&=b\big(1-b\,\alpha(\Gt_{zbz^*}(c))\big)^{-1} \label{eq:Gzcz1} \\
&=b\big(1-b\,\alpha(c(1-c\,\beta(\Gt_{z^*cz}(b)))^{-1})\big)^{-1}. \label{eq:Gzcz2}
\end{align}
Moreover, the $B$--valued R--transform of $z^*cz$ is given by
\begin{equation}\label{eq:Rzcz}
R_{z^*cz}(b)=\alpha(c(1-c\,\beta(b))^{-1}).
\end{equation}
\end{prop}
\begin{proof}
Using cummulants to evaluate $*$--moments of $z$ as in Remark~\ref{rmk:mom}, we have
for $n\ge1$,
\[
E(b(z^*czb)^n)=\sum_{\NC_2(2n)}b\big(\pi\{z^*c,zb,\ldots,z^*c,zb\}\big),
\]
Any $\pi\in\NC_2(2n)$ can be uniquely written as
\begin{equation}\label{eq:pi2k}
\pi=\{\{1,2k\}\}\cup\pi'\cup\pi''
\end{equation}
for some $k\in\{1,\ldots,n\}$,
$\pi'\in\NC_2(\{2,\ldots,2k-1\})$ and $\pi''\in\NC_2(\{2k+1,\ldots,2n\})$,
where $\NC_2(S)$ for a subset $S\subseteq\Ints$ is the set of all non--crossing
pair partitions of $S$, and where we set $\NC_2(\emptyset)=\{\emptyset\}$.
Moreover, for any $k\in\{1,\ldots,n\}$, the map
\[
(\pi',\pi'')\mapsto\{\{1,2k\}\}\cup\pi'\cup\pi''
\]
is a bijection from $\NC_2(\{2,\ldots,2k-1\})\times\NC_2(\{2k+1,\ldots,2n\})$
to $\{\pi\in\NC_2(2n)\mid\{1,2k\}\in\pi\}$.
Finally, with $\pi$ as in~\eqref{eq:pi2k}, we have
\begin{equation*}
\pi\{z^*c,zb,\ldots,z^*c,zb\}
=\alpha(c(\pit'\{zb,z^*c,\ldots,zb,z^*c\}))b(\pit''\{z^*c,zb,\ldots,z^*c,zb\}),
\end{equation*}
where $\pit'\in\NC_2(2k-2)$ and $\pit''\in\NC_2(2n-2k)$ are obtained from $\pi'$ and $\pi''$
by left shifting by $1$ and, respectively, $2k$.
Therefore,
\begin{align*}
E(b(z^*czb)^n)
&=\sum_{k=1}^n\sum_{\substack{\pit'\in\NC_2(2k-2)\\ \pit''\in\NC_2(2n-2k)}}
 \begin{aligned}[t]
 b\,\alpha(c(\pit'\{&zb,z^*c,\ldots,zb.z^*c\})) \\ &b\,(\pit''\{z^*c,zb,\ldots,z^*c,zb\})
 \end{aligned} \\
&=\sum_{k=1}^nb\alpha(E(c(zbz^*c)^{k-1}))E(b(z^*czb)^{n-k}).
\end{align*}
So we have
\begin{align*}
\Gt_{z^*cz}(b)&=\sum_{n=0}^\infty E(b(z^*czb)^n)
=b+b\sum_{n=1}^\infty\sum_{k=1}^n\alpha(E(c(zbz^*c)^{k-1}))E(b(z^*czb)^{n-k}) \\
&=b+b\alpha\bigg(\sum_{r=0}^\infty E(c(zbz^*c)^r)\bigg)\bigg(\sum_{s=0}^\infty E(b(z^*czb)^s)\bigg) \\
&=b+b\alpha(\Gt_{zbz^*}(c))\Gt_{z^*cz}(b).
\end{align*}
Solving yields~\eqref{eq:Gzcz1}.
Since $z^*$ is $B$--circular with covariance $(\beta,\alpha)$, we have
\[
\Gt_{zbz^*}(c)=c\big(1-c\beta(\Gt_{z^*cz}(b))\big)^{-1}.
\]
Plugging this into~\eqref{eq:Gzcz1} yields~\eqref{eq:Gzcz2}.

By~\cite[Thm.\ 4.1.12]{Sp}, which is due to Voiculescu~\cite{V} in a slightly different guise,
the $B$--valued R--transform is
\[
R_{z^*cz}(b)+b^{-1}=K(b)^{-1},
\]
where
\[
\Gt_{z^*cz}(K(b))=K(\Gt_{z^*cz}(b))=b.
\]
Therefore,
\[
b=\Gt_{z^*cz}(K(b))=K(b)\big(1-K(b)\alpha(c(1-c\beta(b))^{-1}))^{-1}.
\]
Solving yields
\[
K(b)^{-1}=b^{-1}+\alpha(c(1-c\beta(b))^{-1}),
\]
and this gives immediately~\eqref{eq:Rzcz}.
\end{proof}

\begin{remark}
The formula~\eqref{eq:Gzcz2} gives a continued--fraction--type expansion:
\[
\Gt_{z^*cz}(b)=\cfrac b{1-b\alpha\bigg(\cfrac c{1-c\beta\bigg(\cfrac b{1-b\alpha\Big(\cfrac c{1-\cdots}\Big)}\bigg)}\bigg)}.
\]
When $c=1$ and when $b=\zeta^{-1}\in\Cpx$, we therefore have
\[
G_{z^*z}(\zeta)=\Gt_{z^*z}(\zeta^{-1})
=\cfrac1{\zeta-\alpha\bigg(\cfrac1{\zeta-\beta\bigg(\cfrac1{\zeta-\alpha\Big(\cfrac1{\zeta-\cdots}\Big)}\bigg)}\bigg)}.
\]
\end{remark}

\begin{prop}\label{prop:pos}
Let $B$ be a unital $C^*$--algebra and let $\alpha,\beta:B\to B$ be
$\Cpx$--linear maps.
Then a $B$--circular operator with covariance $(\alpha,\beta)$
can be realized in a $B$--valued $C^*$-- noncommutative
probability space $(A,E)$ if and only if $\alpha$ and $\beta$ are completely positive.
\end{prop}
\begin{proof}
Necessity follows from the complete positivity of a projection $E:A\to B$ onto a C$^*$--subalgebra, which was proved by Tomiyama~\cite{Tom}.
Sufficiency follows from results of Speicher~\cite{Sp}.
Indeed, using Proposition~\ref{prop:Bsemicirc}, complete positivity of $\alpha$ and $\beta$ implies that
the covariance matrix $\eta:B\to M_2(B)=B\otimes M_2(\Cpx)$ given by
\[
\eta(b)=\frac12\left(\begin{matrix}
\alpha(b)+\beta(b)&i(\beta(b)-\alpha(b))\\
i(\alpha(b)-\beta(b))&\alpha(b)+\beta(b)
\end{matrix}\right)
=\alpha(b)\otimes p
+\beta(b)\otimes(1-p),
\]
where $p=\frac12\left(\begin{smallmatrix}
1&-i\\
i&1
\end{smallmatrix}\right)\in M_2(\Cpx)$, is completely positive;
by~\cite[Thm.\ 4.3.1]{Sp}, the restriction of $E$ to the $*$--algebra $\Afr$
generated by $\{z\}\cup B$ is positive and
by~\cite[Rmk.\ 4.3.2]{Sp}, the $B$--Gaussian random variables with covariance matrix $\eta$ can be
realized in a $B$--valued C$^*$--noncommutative probability space.
\end{proof}

The following exactness result is a direct consequence of Proposition~\ref{prop:Bsemicirc}, \cite[Cor.\ 2.3]{DS}
and the fact that exactness passes to C$^*$--subalgebras.
\begin{prop}\label{prop:exact}
Let $B$ be an exact C$^*$--algebra, let $(A,E)$ be a $B$--valued C$^*$--noncommutative probability
space and let $z\in A$ be a $B$--circular element.
Then the C$^*$--algebra $C^*(B\cup\{z\})$ is exact.
\end{prop}

\begin{lemma}\label{lem:normx}
Let $B$ be a unital C$^*$--algebra and suppose $(A,E)$ is a $B$--valued C$^*$--noncommutative probability space
with $E$ faithful.
For $x\in A$, we have
\[
\|x\|=\limsup_{n\to\infty}\|E((x^*x)^n)\|^{1/2n}.
\]
\end{lemma}
\begin{proof}
We may without loss of generality assume $A$ and $B$ are separable.
Let $\phi$ be a faithful state on $B$.
Then $\phi\circ E$ is a faithful state on $A$.
Hence
\[
\|x\|=\|x^*x\|^{1/2}=\lim_{n\to\infty}(\phi\circ E((x^*x)^n))^{1/2n}
\le\limsup_{n\to\infty}\|E((x^*x)^n)\|^{1/2n}\le\|x\|.
\]
\end{proof}

\begin{prop}\label{prop:normz}
Let $B$ be a unital C$^*$--algebra and suppose $(A,E)$ is a $B$--valued C$^*$--noncommutative probability space
with $E$ faithful.
Let $\alpha$ and $\beta$ be completely positive maps from $B$ to itself.
Suppose $z\in A$ is a $B$--circular element with covariance $(\alpha,\beta)$.
Then
\begin{equation}\label{eq:normz}
\max(\|\alpha\|,\|\beta\|)^{1/2}\le\|z\|\le2\max(\|\alpha\|,\|\beta\|)^{1/2}.
\end{equation}
\end{prop}
\begin{proof}
Let $K=\max(\|\alpha\|,\|\beta\|)$.
Using the recursive formula~\eqref{eq:piinduc} for evaluating the bracketing~\eqref{eq:pibrak},
one sees by induction on $n\ge1$ that
\[
\|\pi\{z^{s(1)}b_1,\ldots,z^{s(2n)}b_{2n}\}\|\le K^{n}\max(\|b_1\|,\ldots,\|b_{2n}\|).
\]
From~\eqref{eq:zmom}, we therefore get
\[
\|E(z^{s(1)}\cdots z^{s(n)})\|\le K^n(\#\NC_2(2n))
\]
whenever $s(1),\ldots,s(2n)\in\{1,*\}$, where $\#\NC_2(n)$ is the number of
non-crossing pair partitions of $\{1,\ldots,n\}$.
Therefore, $\|E((z^*z)^n)\|\le K^n\frac1{n+1}\binom{2n}n$.
Now Lemma~\ref{lem:normx} and the asymptotics of Catalan numbers yield the upper bound in~\eqref{eq:normz}.

For the lower bound, we have
\begin{align*}
\|z^*z\|&\ge\|E(z^*z)\|=\|\alpha(1)\|=\|\alpha\| \\
\|zz^*\|&\ge\|E(zz^*)\|=\|\beta(1)\|=\|\beta\|. \\
\end{align*}
\end{proof}

\section{Hyperinvariant subspaces for certain $L^\infty([0,1])$--circular operators}
\label{sec:hisp}

For completeness, we provide a proof of the following well known characterization
of normal, completely positive maps from $L^\infty(X,\mu)$ to itself, for $\mu$ a probability measure.
This may be compared to~\cite[Ex.\ 2.8]{Sh}, where, however, some conditions are different.
\begin{lemma}\label{lem:cpLinf}
Let $\mu$ be a probability measure on a measurable space $(X,\MEu)$.
Let $\pi_2:X\times X\to X$ be the coordinate projection $\pi_2(x,y)=y$.
Let $\eta$ be a finite, positive measure on $(X\times X,\,\MEu\otimes\MEu)$ and assume that
the push--forward measure ${\pi_2}_*\eta$ is absolutely continuous with respect to $\mu$
and that the Radon--Nikodym derivative $\frac{d({\pi_2}_*\eta)}{d\mu}$ is bounded.
Then there is a (unique) normal, completely positive map 
\begin{equation*}
\alpha_\eta:L^\infty(X,\mu)\to L^\infty(X,\mu)
\end{equation*}
such that for all $h\in L^1(X,\mu)$,
\begin{equation}\label{eq:alphaf}
\int_X(\alpha_\eta f)(y)h(y)d\mu(y)=\int_{X\times X}f(x)h(y)d\eta(x,y).
\end{equation}
We may formally write
\[
(\alpha_\eta f)(y)=\int_Xf(x)\eta(dx,y).
\]
Conversely, every normal, completely positive map
\begin{equation}\label{eq:alpha}
\alpha:L^\infty(X,\mu)\to L^\infty(X,\mu)
\end{equation}
arises in this way from a measure $\eta$, and
\begin{equation}\label{eq:RN}
\frac{d({\pi_2}_*\eta)}{d\mu}=\alpha(1_X).
\end{equation}
Consequently,
\begin{equation}\label{eq:alphanm}
\|\alpha\|=\|\alpha(1_X)\|_\infty=\bigg\|\frac{d({\pi_2}_*\eta)}{d\mu}\bigg\|_\infty.
\end{equation}
\end{lemma}
\begin{proof}
We have
\begin{gather*}
\int|f(x)h(y)|d\eta(x,y)\le\|f\|_\infty\int|h(y)|d\eta(x,y)
=\|f\|_\infty\int|h(y)|d({\pi_2}_*\eta)(y) \\
\le\|f\|_\infty\bigg\|\frac{d({\pi_2}_*\eta)}{d\mu}\bigg\|_\infty\int|h(y)|d\mu(y)
=\|f\|_\infty\bigg\|\frac{d({\pi_2}_*\eta)}{d\mu}\bigg\|_\infty\|h\|_{L^1(\mu)}.
\end{gather*}
So~\eqref{eq:alphaf} uniquely defines an element $\alpha_\eta f$ of $L^\infty(X,\mu)$.
Cleary the map $\alpha_\eta$ is positive (therefore, completely positive) and normal.

Conversely, given a normal, completely positive map $\alpha$ as in~\eqref{eq:alpha},
for $E_1,E_2\in\MEu$ define
\begin{equation}\label{eq:etaE}
\eta(E_1\times E_2)=\int_{E_2}(\alpha(1_{E_1}))(y)d\mu(y).
\end{equation}
Using positivity and normality, $\eta$ is seen to extend to a finite, positive measure
on $X\times X$.
Finally, ${\pi_2}_*\eta(E)=\eta(X\times E)$ and from~\eqref{eq:etaE} we get that ${\pi_2}_*\eta$ is $\mu$--absolutely continuous
and~\eqref{eq:RN} holds.
Now equation~\eqref{eq:alphanm} follows directly.
\end{proof}

If we desire a completely positive map $\beta_\eta:L^\infty(X,\mu)\to L^\infty(X,\mu)$ such that
\[
\tau(\alpha_\eta(f)g)=\tau(f\beta_\eta(g))
\]
for all $f,g\in L^\infty(X,\mu)$, where $\tau(\cdot)=\int\cdot\,d\mu$, (cf.\ Proposition~\ref{prop:tracial}),
then $\beta_\eta$ will need to satisfy
\begin{equation*}
\int_X(\beta_\eta g)(x)h(x)d\mu(x)=\int_{X\times X}h(x)g(y)d\eta(x,y)
\end{equation*}
for all $h\in L^1(X,\mu)$, and we will need also ${\pi_1}_*\eta$ to be absolutely continuous with respect to
$\mu$ and have bounded Radon--Nikodym derivative, where $\pi_1:X\times X\to X$ is the other coordinate projection.
We may formally write
\[
(\beta_\eta g)(x)=\int_Xg(y)\eta(x,dy).
\]

Consider $\Dc=L^\infty([0,1])$ with trace $\tau$ given by integration with respect to Lebesgue
measure.
We will study $\Dc$--circular operators in a $\Dc$--valued W$^*$--noncommutative probability space
$(\Mcal,\Ec)$ such that $\tau\circ\Ec$ is a normal, faithful, tracial state on $\Mcal$.
In light of the above discussion, this class of operators is precisely the class delineated
below.

\begin{defi}\label{def:zeta}
If $\eta$ is a finite, Borel measure on $[0,1]^2$ whose push--forward
measures ${\pi_i}_*\eta$ under both coordinate projections $\pi_1,\pi_2:[0,1]^2\to[0,1]$
are absolutely continuous with respect to Lebesgue measure and have bounded
Radon--Nikodym derivative, let $z_\eta$ be a $\Dc$--circular operator with
covariance $(\alpha_\eta,\beta_\eta)$, where, for all $h\in L^1([0,1])$,
\begin{equation}
\label{eq:alphabetaeta}
\begin{aligned}
\int_0^1(\alpha_\eta f)(y)h(y)dy&=\int_{[0,1]^2}f(x)h(y)d\eta(x,y), \\
\int_0^1(\beta_\eta g)(x)h(x)dx&=\int_{[0,1]^2}h(x)g(y)d\eta(x,y).
\end{aligned}
\end{equation}
\end{defi}
It follows from Propositions~\ref{prop:pos} and~\ref{prop:tracial} that such a $\Dc$--circular
operator $z_\eta$ exists in a tracial $\Dc$--valued C$^*$--noncommutative probability space, and the
Gelfand--Naimark--Segal construction then yields $z_\eta$ in
a $\Dc$--valued W$^*$--noncommutative probability space
$(\Mcal,\Ec)$, with $\tau\circ\Ec$ a normal, faithful, tracial state on $\Mcal$.
We may also use $\tau$ to denote the faithful trace $\tau\circ\Ec$ on $\Mcal$, and for $a\in\Mcal$,
we let $\|a\|_2=\tau(a^*a)^{1/2}$, as usual.
\begin{lemma}\label{lem:2norm}
$\|z_\eta\|_2=\eta([0,1]^2)^{1/2}$.
\end{lemma}
\begin{proof}
We have
\[
\|z_\eta\|_2^2=\tau(z_\eta^*z_\eta)=\tau\circ\Ec(z_\eta^*z_\eta)=\tau(\alpha_\eta(1))
=\int_0^1(\alpha_\eta(1))(y)dy=\eta([0,1]^2),
\]
where the last equality is from~\eqref{eq:alphabetaeta}.
\end{proof}

\begin{notation}\label{not:chi}
For Borel subsets $A$ of $\Reals$, we will use $1_A$ for the characteristic function
of $A$, for example as in $1_A(S)$, when applied via the Borel functional calculus
to a self--adjoint operator $S\in\Mcal$.
On the other hand, the notation $\chi_A$ for $A\subseteq[0,1]$ will be used for the characteristic
function of $A$ considered as an element of $\Dc=L^\infty([0,1])\subseteq\Mcal$.
\end{notation}

\begin{lemma}\label{lem:chiAB}
If $A$ and $B$ are Borel subsets of $[0,1]$ and if $\eta(A\times B)=0$, then
$\chi_Az_\eta\chi_B=0$.
\end{lemma}
\begin{proof}
We have
\begin{align*}
\tau\circ\Ec((\chi_Az_\eta\chi_B)^*(\chi_Az_\eta\chi_B))&=\tau(\chi_B\Ec(z_\eta^*\chi_Az_\eta))
=\tau(\chi_B\alpha_\eta(\chi_A)) \\
&=\int_{[0,1]^2}\chi_A(x)\chi_B(y)d\eta(x,y)=\eta(A\times B)=0.
\end{align*}
\end{proof}

\begin{lemma}\label{lem:Ez}
Suppose for some $0\le a<1$, the restriction of $\eta$ to $[a,1]\times[0,1]$
is absolutely continuous with respect to Lebesgue measure,
has bounded Radon---Nikodym derivative and is supported in
$\{(s,t)\mid a\le s\le t\le1\}$.
Then there is $K>0$ such that for all $n\in\Nats$ and all $\mu\in[0,\frac{K^n(1-a)^n}{n!}]$, we have
\begin{equation*}
1_{[0,\mu]}(\Ec(z_\eta^n(z_\eta^*)^n))\ge\chi_{[\rho,1]},
\end{equation*}
where 
\begin{equation}\label{eq:rho}
\rho=1-\frac{(n!\,\mu)^{1/n}}K.
\end{equation}
Consequently, if we fix $\gamma\in[a,1)$, then letting
$\zeta=(\chi_{[\gamma,1]})\widehat{\;}\in L^2(\Dc,\tau\restrict_\Dc)$
and letting
\begin{equation}\label{eq:mun}
\mu_n=\frac{K^n(1-\gamma)^n}{n!},
\end{equation}
we have
\[
1_{[0,\mu_n]}(\Ec(z_\eta^n(z_\eta^*)^n))\zeta=\zeta
\]
for all $n\in\Nats$.
\end{lemma}
\begin{proof}
Let $H$ be the Radon--Nikodym derivative of the restriction of $\eta$
to $[a,1]\times[0,1]$ with respect to Lebesgue measure and let $K>0$
be at least as large as
the essential supremum $\|H\|_\infty$ of $H$.
Since $\eta([a,1]\times[0,a])=0$, from Lemma~\ref{lem:chiAB} we have
$\chi_{[a,1]}z_\eta=\chi_{[a,1]}z_\eta\chi_{[a,1]}$ and
consequently,
\begin{equation}\label{eq:Ecazeta}
\Ec(z_\eta^n(z_\eta^*)^n)\ge\chi_{[a,1]}\Ec(z_\eta^n(z_\eta^*)^n)
=\Ec((\chi_{[a,1]}z_\eta)^n\chi_{[a,1]}(z_\eta^*\chi_{[a,1]})^n).
\end{equation}
Let $f_n$ denote the right--hand--side of~\eqref{eq:Ecazeta},
with $f_0=\chi_{[a,1]}$.
Then by the nested evaluation described in Remark~\ref{rmk:mom}, we get
$f_{n+1}=\chi_{[a,1]}\beta_\eta(f_n)$, for all $n\ge0$.
We therefore have, whenever $a\le x\le1$,
\[
0\le f_{n+1}(x)=(\beta_\eta f_n)(x)=\int_a^1f_n(y)H(x,y)dy\le K\int_x^1f_n(y)dy.
\]
It follows by induction on $n\ge0$ that
\[
0\le f_n(x)\le K^n\frac{(1-x)^n}{n!},\qquad(a\le x\le 1).
\]
If $0\le\mu\le\frac{K^n(1-a)^n}{n!}$, then 
\[
1_{[0,\mu]}(\Ec(z_\eta^n(z_\eta^*)^n))\ge1_{[0,\mu]}(f_n)
\ge1_{[0,\mu]}(\frac{K^n(1-x)^n}{n!})=\chi_{[\rho,1]},
\]
where $\frac{K^n(1-\rho)^n}{n!}=\mu$, i.e.\ where~\eqref{eq:rho} holds.
The remaining assertions follow directly.
\end{proof}

\begin{lemma}\label{lem:zcd}
Let $0\le c<d\le1$ and suppose
\begin{equation}
\label{eq:etacd}
\eta([d,1]\times[c,d])=0=\eta([c,1]\times[0,c]).
\end{equation}
Let $\phi:[c,d]\to[0,1]$ be $\phi(x)=(x-c)/(d-c)$.
Let
\[
\eta'=(d-c)^{-1}(\phi\times\phi)_*(\eta\restrict_{[c,d]^2})
\]
be the measure on $[0,1]^2$ that is $(d-c)^{-1}$ times the push--forward under $\phi\times\phi$
of the restriction of $\eta$ to $[c,d]\times[c,d]$.
Then whenever $0\le\theta\le d-c$ and $n\in\Nats$, we have
\begin{equation}
\label{eq:zeta'}
s_\theta(z_\eta^n)\ge s_{\frac\theta{d-c}}(z_{\eta'}^n).
\end{equation}
\end{lemma}
\begin{proof}
Let $p_1=\chi_{[0,c]}$ and $p_2=\chi_{[c,d]}$.
From~\eqref{eq:etacd} and Lemma~\ref{lem:chiAB}, $z_\eta p_2=(p_1+p_2)z_\eta p_2$ and $z_\eta p_1=p_1z_\eta p_1$, so
$p_2z_\eta^np_2=(p_2z_\eta p_2)^n$.
Consequently, by Lemma~\ref{lem:s},
\begin{equation}
\label{eq:p2zeta}
s_\theta(z_\eta^n)\ge s_{\frac\theta{d-c}}((p_2z_\eta p_2)^n;(d-c)^{-1}(\tau\circ\Ec)\restrict_{p_2\Mcal p_2}).
\end{equation}
By Proposition~\ref{prop:properties}(iv), $p_2z_\eta p_2$ is a $p_2\Dc$--circular element
in $(p_2\Mcal p_2,\Ec\restrict_{p_2\Mcal p_2})$ with covariance $(\alphat,\betat)$, where for $f\in p_2\Dc=L^\infty([c,d])$,
$\alphat(f)=p_2\alpha_\eta(f)$ and $\betat(f)=p_2\beta_\eta(f)$.

Consider the isomorphism $\phit:L^\infty([0,1])\to L^\infty([c,d])$ given by $\phit(f)=f\circ\phi$.
Using this identification, we may regard $(p_2\Mcal p_2,\Ec\restrict_{p_2\Mcal p_2})$ as an $L^\infty([0,1])$--valued
W$^*$--noncommutative probability space and $p_2z_\eta p_2$ as an $L^\infty([0,1])$--circular operator with covariance
$(\phit^{-1}\circ\alphat\circ\phit,\phit^{-1}\circ\betat\circ\phit)$.
Let $f\in L^\infty([0,1])$ and $h\in L^1([0,1])$.
Then
\begin{align*}
\int_0^1((\phit^{-1}\circ\alphat\circ\phit)f)(t)h(t)dt
&=\int_0^1((\alphat\circ\phit)f)(c+(d-c)t)h(t)dt \\
&=(d-c)^{-1}\int_c^d((\alphat\circ\phit)f)(y)h(\frac{y-c}{d-c})dy \\
&=(d-c)^{-1}\int_{[c,d]^2}(\phit f)(x)h(\frac{y-c}{d-c})d\eta(x,y) \\
&=(d-c)^{-1}\int_{[c,d]^2}(f\circ\phi)(x)(h\circ\phi)(y)d\eta(x,y) \\
&=(d-c)^{-1}\int_{[0,1]^2}f(s)h(t)d((\phi\times\phi)_*\eta)(s,t) \\
&=\int_0^1(\alpha_{\eta'}f)(t)h(t)dt.
\end{align*}
Thus, $\phit^{-1}\circ\alphat\circ\phit=\alpha_{\eta'}$.
Similarly, we have $\phit^{-1}\circ\betat\circ\phit=\beta_{\eta'}$.
Hence, $p_2z_\eta p_2$ is identified with $z_{\eta'}$.
Finally, note that the trace $(d-c)^{-1}\tau\circ\Ec\restrict_{p_2\Mcal p_2}$ is equal to
$\tau\circ\phit^{-1}\circ\Ec\restrict_{p_2\Mcal p_2}$, and~\eqref{eq:zeta'} follows from~\eqref{eq:p2zeta}.
\end{proof}

\begin{thm}\label{thm:zhis}
Consider an $L^\infty([0,1])$--circular operator $z_\eta$ as described in Definition~\ref{def:zeta}.
Suppose\renewcommand{\labelenumi}{(\roman{enumi})}
\begin{enumerate}
\item for some $0\le a<1$, the restriction of $\eta$ to $\{(s,t)\mid a\le s\le t\le 1\}$ is
less than or equal to $R$ times Lebesgue measure, for some $R<\infty$;
\item for some $0\le c<d\le1$, the restriction of $\eta$ to $\{(s,t)\mid c\le s\le t\le d\}$ is
$r$ times Lebesgue measure for some $r>0$;
\item $\eta$ vanishes on
\begin{align*}
([c,1]\times[0,c])&\cup([d,1]\times[c,d])\cup([a,1]\times[0,1]) \\
&\cup\{(s,t)\mid c\le t\le s\le d\}\cup\{(s,t)\mid a\le t\le s\le a\}.
\end{align*}
\end{enumerate}
Then $z_\eta$ has a nontrivial, hyperinvariant subspace.
\end{thm}
\begin{proof}
Note that we may without loss of generality take $d<a$.
The conditions of the theorem are illustrated in Figure~\ref{fig:etagencond}.
\begin{figure}[b]
\begin{picture}(200,200)(0,-30)

  \linethickness{0.7pt}
  \drawline(0,0)(0,160)(160,160)(160,0)(0,0) 

  \linethickness{0.5pt}

  \drawline(0,135)(75,135)(75,52) 
  \drawline(75,85)(25,135)
  \put(56,116){$r$}

  \drawline(75,52)(160,52) 
  \drawline(160,0)(108,52)
  \put(130,35){$\le R$}

  \put(25,100){$0$}
  \put(45,40){$0$}
  \put(100,15){$0$}

  \put(10,145){$*$}
  \put(90,70){$*$}
  \put(120,116){$*$}

  \drawline(160,160)(190,160) 
  \drawline(187,162)(190,160)(187,158) 
  \put(175,165){$t$}
  \drawline(160,160)(160,165) 
  \put(154,168){$1$}
  \drawline(108,160)(108,165) 
  \put(102,168){$a$}
  \drawline(75,160)(75,165) 
  \put(69,168){$d$}
  \drawline(25,160)(25,165) 
  \put(19,168){$c$}

  \drawline(0,0)(0,-30) 
  \drawline(-2,-27)(0,-30)(2,-27) 
  \put(-15,-20){$s$}
  \drawline(0,0)(-5,0) 
  \put(-13,3){$1$}
  \drawline(0,52)(-5,52) 
  \put(-13,55){$a$}
  \drawline(0,85)(-5,85) 
  \put(-13,88){$d$}
  \drawline(0,135)(-5,135) 
  \put(-13,138){$c$}

\end{picture}
\caption{Conditions on $\eta$ from Theorem~\ref{thm:zhis}.}\label{fig:etagencond}
\end{figure}
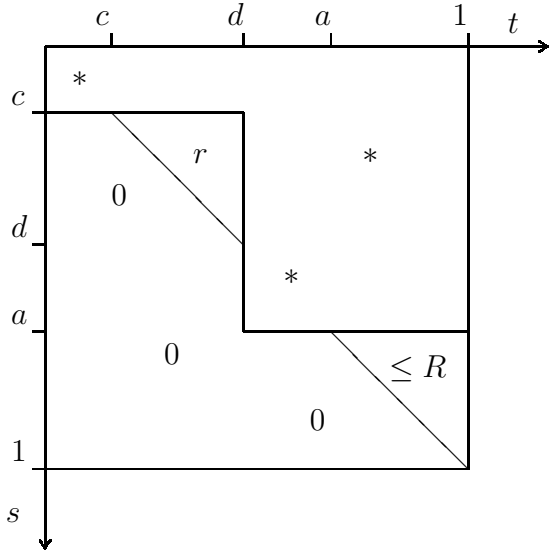

By Lemma~\ref{lem:zcd}, if $0<\theta<d-c$ and $n\in\Nats$, then $s_\theta(z_\eta^n)\ge s_{\frac\theta{d-c}}(z_{\eta'}^n)$,
where $\eta'$ is $r(d-c)$ times the Lebesgue measure supported on $\{(x,y)\mid 0\le x\le y\le 1\}$.
Therefore, cf.\ Examples~\ref{ex:Bcirc}(iii), $z_{\eta'}$ is a $\DT(\delta_0,\sqrt{r(d-c)})$--operator, i.e.\ is
$\sqrt{r(d-c)}$ times a $\DT(\delta_0,1)$--operator $T$.
By \'Sniady's result~\cite{Sn} on $*$--moments of $T$, it follows that $(T^*)^nT^n$ and $(\frac1nT^*T)^n$ have the same
$*$--moments.
Hence, for any $n\in\Nats$ and $0<\sigma<1$,
\[
s_\sigma(T^n)=s_\sigma((T^*)^nT^n)^{1/2}=n^{-n/2}s_\sigma((T^*T)^n)^{1/2}=n^{-n/2}s_\sigma(T^*T)^{n/2}=n^{-n/2}s_\sigma(T)^n.
\]
Hence,
\[
s_\theta(z_\eta^n)\ge s_{\frac\theta{d-c}}(z_{\eta'}^n)=\bigg(\frac{r(d-c)}n\bigg)^{n/2}s_{\frac\theta{d-c}}(T)^n.
\]
By~\cite{DH2}, the operator $T$ has trivial kernel
(in fact, the distribution of $T^*T$ was explicitly determined there).
Fixing any $\theta\in(0,d-c)$, we get $s_{\frac\theta{d-c}}(T)\ne0$, and
\[
s_\theta(z_\eta^n)\ge(\frac\alpha n)^{n/2}
\]
for some $\alpha>0$.

We may apply Lemma~\ref{lem:Ez} to $z_\eta$.
Let $K$ be as in that lemma, and choose $\gamma$ sufficiently close to $1$ so that $K(1-\gamma)\le\frac\alpha e$.
Then choosing $\mu_n$ as in~\eqref{eq:mun} and using Stirling's formula for $n!$, we have
\[
\limsup_{n\to\infty}\frac{\mu_n}{s_\theta(z_\eta^n)^2}\le\limsup_{n\to\infty}\frac{n^n}{n!}\bigg(\frac{K(1-\gamma)}\alpha\bigg)^n
=\limsup_{n\to\infty}\frac1{c_n\sqrt n}\bigg(\frac{eK(1-\gamma)}\alpha\bigg)^n=0,
\]
where $c_n$ converges to a strictly postive number.
Therefore, Theorem~\ref{thm:FJKP} applies, with $p=0$, and yields a nontrivial hyperinvariant subspace for $z_\eta$.
\end{proof}

\section{$L^\infty([0,1])$--circular operators in free group factors}
\label{sec:inLF3}

In this section, we construct an $L^\infty([0,1])$--circular operator $z_\eta$, as in Definition~\ref{def:zeta},
inside of a free group factor, when $\eta$ is assumed to be absolutely continuous with respect to
Lebesgue measure.
This construction parallels what was done in~\cite[\S4]{DH2} for the quasinilpotent DT--operator.

As in the previous section, $\Dc$ will denote $L^\infty([0,1])$ with trace $\tau$ given by
integration with respect to Lebesgue measure, $(\Mcal,\Ec)$ will be a $\Dc$--valued
W$^*$--noncommutative probability space such that $\tau\circ\Ec$ is a normal, faithful,
tracial state on $\Mcal$,
and $z_\eta$ will be a $\Dc$--circular operator in $(\Mcal,\Ec)$ with covariance $(\alpha_\eta,\beta_\eta)$,
described by equations~\eqref{eq:alphabetaeta}.
We write $\tau$ also for the trace $\tau\circ\Ec$.
Let $H\in L^1([0,1])$ be the Radon--Nikodym derivative of $\eta$ with respect to Lebesgue measure.
Then
\begin{equation}\label{eq:alphabetaH}
\begin{aligned}
(\alpha_\eta f)(t)&=\int_0^1H(s,t)f(s)ds \\
(\beta_\eta f)(t)&=\int_0^1H(t,u)f(u)du.
\end{aligned}
\end{equation}
Moreover, the push--forward measures $(\pi_i)_*\eta$ of $\eta$ under the coordinate projections $\pi_1$ and $\pi_2$
are absolutely continuous with respect to Lebesgue measure and have Radon--Nikodym derivatives
equal to the {\em coordinate expectations} $CE_1(H)$ and $CE_2(H)$,
respectively, given by
\begin{equation}\label{eq:h1h2}
\begin{aligned}
CE_1(H)(x)&=\int_0^1H(x,y)dy \\
CE_2(H)(y)&=\int_0^1H(x,y)dx.
\end{aligned}
\end{equation}
Thus, the assumption on $\eta$ from Definition~\ref{def:zeta} is that $CE_1(H)$ and $CE_2(H)$ are essentially
bounded, and then from Proposition~\ref{prop:normz} and equation~\eqref{eq:alphanm} of
Lemma~\ref{lem:cpLinf}, we have
\begin{equation}\label{eq:CEHnorm}
\|z_\eta\|\le2\max(\|CE_1(H)\|_\infty,\|CE_2(H)\|_\infty)^{1/2}.
\end{equation}

\begin{defi}\label{def:Mblock}
Let $w\in L^\infty([0,1]^2)$.
We will say $w$ is in {\em regular block form} if $w$ is constant on all blocks in the regular $n\times n$ lattice
superimposed on $[0,1]^2$, for some $n$, i.e.\ if there are $n\in\Nats$ and $w_{ij}\in\Cpx$, ($1\le i,j\le n$),
such that $w(s,t)=w_{ij}$ whenever
$\frac{i-1}n\le s<\frac in$ and $\frac{j-1}n\le t<\frac jn$, for all integers $1\le i,j\le n$.
For specificity, we may then say that $w$ is in $n\times n$ regular block form.
Let $a\in\Mcal$.
Then we set
\[
M(w,a)=\sum_{i,j=1}^nw_{ij}p_iap_j\in\Mcal,
\]
where $p_i=\chi_{[\frac{i-1}n,\frac in]}\in\Dc$.
\end{defi}

The following properties are straightforward.
\begin{lemma}\label{lem:Mblock}
\renewcommand{\labelenumi}{(\alph{enumi})}
\begin{enumerate}
\item $M(w,a)$ is independent of the choice of $n$ so long as $w$ is in $n\times n$ regular block form.
\item $M(w,\zeta a_1+a_2)=\zeta M(w,a_1)+M(w,a_2)$ for $a_1,a_2\in\Mcal$ and $\zeta\in\Cpx$.
\item If $w^{(1)},\,w^{(2)}\in L^\infty([0,1]^2)$ are both in regular block form, then there is $n\in\Nats$
such that both are in $n\times n$ regular block form;
for $\zeta\in\Cpx$, we then have
\[
M(\zeta w^{(1)}+w^{(2)},a)=\zeta M(w^{(1)},a)+M(w^{(2)},a).
\]
\end{enumerate}
\end{lemma}

For the rest of this section, we will suppose that
$z\in\Mcal$ is (scalar) circular with respect to $\tau$ and satisfies $\tau(z)=0$ and $\tau(z^*z)=1$
and that $\Dc$ and $z$ are $*$--free (over $\Cpx$) with respect to $\tau$.
Therefore, $W^*(\Dc\cup\{z\})\cong L(\Fb_3)$.

\begin{lemma}\label{lem:Mblockcirc}
Let $w\in L^\infty([0,1]^2)$ be in regular block form, let $H=|w|^2$ and let $\eta$ be the
Lebesgue--absolutely--continuous measure on $[0,1]^2$
whose Radon--Nikodym derivative is $H$.
Then $M(w,z)$ is $\Dc$--circular in $(\Mcal,\Ec)$,
with covariance $(\alpha_\eta,\beta_\eta)$.
\end{lemma}
\begin{proof}
For brevity, write $a$ for $M(w,z)$.
Suppose $w$ is in $n\times n$ regular block form with $w_{ij}$ as in Definition~\ref{def:Mblock}.
The equalities
\begin{equation}\label{eq:acov}
\begin{gathered}
\Ec(a^*fa)=\alpha_\eta(f) \\
\Ec(afa^*)=\beta_\eta(f) \\
\Ec(afa)=\Ec(a^*fa^*)=\Ec(a)=\Ec(a^*)=0
\end{gathered}
\end{equation}
for $f\in\Dc$, with $\alpha_\eta$ and $\beta_\eta$ as in~\eqref{eq:alphabetaH},
are easily verified using freeness.
For example,
\begin{align*}
\Ec(a^*fa)&=\sum_{i,j,k=1}^n\overline{w_{ji}}w_{jk}\Ec(p_iz^*p_jfp_jzp_k)
 =\sum_{k=1}^np_k\sum_{j=1}^n|w_{jk}|^2\Ec(z^*p_jfz) \\
&=\sum_{k=1}^np_k\sum_{j=1}^n|w_{jk}|^2\tau(p_jf)
 =\sum_{k=1}^np_k\sum_{j=1}^n|w_{jk}|^2\int_{(j-1)/n}^{j/n}f(s)ds
 =\alpha_\eta(f).
\end{align*}
Suppose $z_1,z_2,\ldots\in\Mcal$ are (scalar) circular elements such that $\tau(z_j)=0$ and $\tau(z_j^*z_j)=1$
and the family
$\Dc,\{z_1\},\,\{z_2\},\ldots$ is $*$--free with respect to $\tau$.
Then the family
$\big(\alg(\{z_j,z_j^*\}\cup\Dc)\big)_{j=1}^\infty$ is free over $\Dc$ with respect to $\Ec$.
Hence, by~\cite[Thm.\ 4.2.4]{Sp}, letting
\[
a_k=\frac1{\sqrt k}\big(M(w,z_1)+\cdots+M(w,z_k)\big)
\]
the pairs $(a_k,a_k^*)$ converge in $\Dc$--valued moments with respect to $\Ec$ to $\Dc$--Gaussian elements
with covariance given by~\eqref{eq:acov} as $k\to\infty$.
In other words, $a_k$ converges in $\Dc$--valued $*$--moments to a $\Dc$--circular element with covariance
$(\alpha_\eta,\beta_\eta)$.
However, by Lemma~\ref{lem:Mblock},
\[
a_k=M\big(w,\frac{z_1+\cdots+z_k}{\sqrt k}\big).
\]
But $z':=\frac{z_1+\cdots+z_k}{\sqrt k}$ is a (scalar) circular element with $\tau(z')=0$ and $\tau((z')^*z')=1$ and with
$\Dc$ and $z'$ $*$--free.
Thus, $a_k$ has the same $\Dc$--valued $*$--moments as $a$
itself, so $a$ is $\Dc$--circular with covariance $(\alpha_\eta,\beta_\eta)$.
\end{proof}

\begin{lemma}\label{lem:Happrox}
Let $H\in L^1([0,1])$, $H\ge0$ and assume the coordinate expectations $CE_1(H)$ and $CE_2(H)$ as in~\eqref{eq:CEHnorm}
are essentially bounded.
Let $\eta$ be the Lebesgue--absolutely--continuous measure on $[0,1]^2$ whose Radon--Nikodym derivative is $H$.
Let $w=\sqrt H$.
Suppose there is a sequence $(w^{(n)})_{n=1}^\infty$ in $L^\infty([0,1]^2)$ such that
\renewcommand{\labelenumi}{(\roman{enumi})}
\begin{enumerate}
\item for each $n$, $w^{(n)}$ is in regular block form,
\item $\lim_{n\to\infty}\|w-w^{(n)}\|_{L^2([0,1]^2)}=0$,
\item letting $H^{(n)}=|w^{(n)}|^2$, both $\|CE_1(H^{(n)})\|_\infty$ and $\|CE_2(H^{(n)})\|_\infty$
remain bound\-ed as $n\to\infty$.
\end{enumerate}
Let $a_n=M(w^{(n)},z)$, with $z$ as in Lemma~\ref{lem:Mblockcirc}.
Then $a_n$ converges in strong--operator topology (in the representation of $\Mcal$ on $L^2(\Mcal,\tau)$)
to an element of $\Mcal$ which is a $\Dc$--circular operator with covariance $(\alpha_\eta,\beta_\eta)$.
\end{lemma}
\begin{proof}
By Lemma~\ref{lem:Mblockcirc} and~\eqref{eq:CEHnorm},
\[
\|a_n\|\le2\max(\|CE_1(H^{(n)})\|_\infty,\|CE_2(H^{(n)})\|_\infty)^{1/2},
\]
so $\|a_n\|$ remains bounded as $n\to\infty$.
From Lemma~\ref{lem:Mblock}, we have
\[
a_n-a_m=M(w^{(n)}-w^{(m)},z).
\]
By Lemma~\ref{lem:Mblockcirc}, $a_n-a_m$ is $\Dc$--circular with covariance corresponding to
the measure on $[0,1]^2$ whose Radon--Nikodym derivative is $|w^{(n)}-w^{(m)}|^2$.
Thus, from Lemma~\ref{lem:2norm} we obtain
\[
\|a_n-a_m\|_2=\|w^{(n)}-w^{(m)}\|_{L^2([0,1]^2)}.
\]
Therefore, $a_n$ is Cauchy in $L^2(\Mcal,\tau)$.
Since $\|a_n\|$ remains bounded, it follows that $a_n$ converges in strong--operator topology to
an element $a$ of $\Mcal$.

It follows, too, that the $\Dc$--valued $*$--moments of $a_n$ converge in strong--operator topology
to those of $a$ as $n\to\infty$, in the sense that
\[
\mbox{s.o.t--}\lim_{n\to\infty}\Ec(d_0a_n^{s(1)}d_1\cdots a_n^{s(k)}d_k)=
\Ec(d_0a^{s(1)}d_1\cdots a^{s(k)}d_k)
\]
for all $k\in\Nats$, $d_0,\ldots,d_k\in\Dc$ and $s(1),\ldots,s(k)\in\{1,*\}$.
Therefore, the $\Dc$--valued free cummulants of $a_n$ converge to those of $a$
in strong--operator topology as $n\to\infty$.
Let $\eta_n$ be the Lebesgue absolutely continuous measure on $[0,1]^2$ whose
Radon--Nikodym derivative is $H^{(n)}$.
By Lemma~\ref{lem:Mblockcirc}, $a_n$ is $\Dc$--circular with covariance $(\alpha_{\eta_n},\beta_{\eta_n})$,
and it follows that $a$ is $\Dc$--circular with covariance $(\alpha,\beta)$, where for $f\in\Dc$,
\begin{align*}
\alpha(f)&=\mbox{s.o.t--}\lim_{n\to\infty}\alpha_{\eta_n}(f) \\
\beta(f)&=\mbox{s.o.t--}\lim_{n\to\infty}\beta_{\eta_n}(f).
\end{align*}
From~\eqref{eq:alphabetaH}, we have
\begin{align*}
(\alpha_{\eta_n} f)(t)&=\int_0^1|w^{(n)}(s,t)|^2f(s)ds \\
(\beta_{\eta_n} f)(t)&=\int_0^1|w^{(n)}(t,u)|^2f(u)du.
\end{align*}
But
\begin{align*}
\|\alpha_{\eta_n}(f)-\alpha_\eta(f)&\|_{L^1([0,1])}=
\int_0^1\big|(\alpha_{\eta_n}(f)-\alpha_\eta(f))(t)\big|dt \\
&=\int_0^1\bigg|\int_0^1\big(|w^{(n)}(s,t)|^2-|w(s,t)|^2\big)f(s)ds\bigg|dt \\
&\le\|\big(|w^{(n)}|^2-|w|^2\big)\|_{L^1([0,1]^2)}\|f\|_\infty \\
&\le\|\big(|w^{(n)}|-|w|\big)\|_{L^2([0,1]^2)}\big(\|w^{(n)}\|_{L^2([0,1]^2)}+\|w\|_{L^2([0,1]^2)}\big)\|f\|_\infty,
\end{align*}
and $\|w^{(n)}\|_{L^2([0,1]^2)}$ remains bounded as $n\to\infty$, while
$\|(|w^{(n)}|-|w|)\|_{L^2([0,1]^2)}\le\|w^{(n)}-w\|_{L^2([0,1]^2)}$ tends to zero.
Therefore, $\alpha=\alpha_\eta$.
Similarly, we find $\beta=\beta_\eta$.
\end{proof}

\begin{thm}\label{thm:inLF3}
Let $H\in L^1([0,1]^2)$ have essentially bounded coordinate expectations $CE_1(H)$ and $CE_2(H)$.
Let $\eta$ be the Lebesgue absolutely continuous measure on $[0,1]^2$ whose Radon--Nikodym derivative
is $H$.
Then there is a $\Dc$--circular operator $z_\eta$ with covariance $(\alpha_\eta,\beta_\eta)$ in $W^*(\Dc\cup\{z\})\cong L(\Fb_3)$.
\end{thm}
\begin{proof}
Let $w=\sqrt H$.
By Lemma~\ref{lem:Happrox}, it will suffice to find a sequence $(w^{(n)})_{n=1}^\infty$ satisfying
hypotheses (i)--(iii) listed there.
For integers $1\le i,j\le n$, let
\[
w_{ij}^{(n)}=n^2\int_{(i-1)/n}^{i/n}\int_{(j-1)/n}^{j/n}w(x,y)dydx
\]
and let $w^{(n)}(s,t)=w_{ij}^{(n)}$ whenever $(s,t)\in[\frac{i-1}n,\frac in)\times[\frac{j-1}n,\frac jn)$.
Then $w^{(n)}$ is in $n\times n$ regular block form; i.e.\ (i) holds.
Let us show 
\begin{equation}\label{eq:wnw}
\lim_{n\to\infty}\|w^{(n)}-w\|_{L^2([0,1]^2)}=0.
\end{equation}
Let $\eps>0$.
There is a continuous function $\wt:[0,1]^2\to[0,\infty)$ such that $\|w-\wt\|_{L^2}<\eps$.
Let
\[
\wt_{ij}^{(n)}=n^2\int_{(i-1)/n}^{i/n}\int_{(j-1)/n}^{j/n}\wt(x,y)dydx
\]
and let $\wt^{(n)}(s,t)=\wt_{ij}^{(n)}$ whenever $(s,t)\in[\frac{i-1}n,\frac in)\times[\frac{j-1}n,\frac jn)$.
Then
\[
\|w^{(n)}-\wt^{(n)}\|_{L^2}^2=n^{-2}\sum_{i,j=1}^n|w_{ij}^{(n)}-\wt_{ij}^{(n)}|^2.
\]
But
\begin{align*}
|w_{ij}^{(n)}-\wt_{ij}^{(n)}|
&\le\int_{(i-1)/n}^{i/n}\int_{(j-1)/n}^{j/n}|w(x,y)-\wt(x,y)|(n^2)dydx \\
&\le\bigg(\int_{(i-1)/n}^{i/n}\int_{(j-1)/n}^{j/n}|w(x,y)-\wt(x,y)|^2(n^2)dydx\bigg)^{1/2} \\
&=n\bigg(\int_{(i-1)/n}^{i/n}\int_{(j-1)/n}^{j/n}|w(x,y)-\wt(x,y)|^2dydx\bigg)^{1/2},
\end{align*}
where the second inequality is because $(n^2)dydx$ is a probability measure on
$[\frac{i-1}n,\frac in)\times[\frac{j-1}n,\frac jn)$.
Therefore,
\[
\|w^{(n)}-\wt^{(n)}\|_{L^2}^2\le\sum_{i,j=1}^n\int_{(i-1)/n}^{i/n}\int_{(j-1)/n}^{j/n}|w(x,y)-\wt(x,y)|^2dydx
=\|w-\wt\|_{L^2}^2<\eps^2.
\]
By uniform continuity of $\wt$, $\lim_{n\to\infty}\|\wt-\wt^{(n)}\|_{L^2}=0$, and using the triangle inequality,
we get $\|w-w^{(n)}\|_{L^2}<3\eps$ for $n$ sufficiently large.
This proves~\eqref{eq:wnw}, namely that hypothesis~(ii) holds.

Finally, for~(iii), letting $H^{(n)}=|w^{(n)}|^2$, we wish to show that $\|CE_1(H^{(n)})\|_\infty$
and $\|CE_2(H^{(n)})\|_\infty$ remain bounded as $n\to\infty$.
We have, for $x\in[\frac{i-1}n,\frac in)$,
\[
(CE_1(H^{(n)}))(x)=\int_0^1|w^{(n)}(x,y)|^2dy=\frac1n\sum_{j=1}^n|w_{ij}^{(n)}|^2.
\]
But
\begin{align*}
|w_{ij}^{(n)}|&=\int_{(i-1)/n}^{i/n}\int_{(j-1)/n}^{j/n}w(x,y)(n^2)dydx \\
&\le\bigg(\int_{(i-1)/n}^{i/n}\int_{(j-1)/n}^{j/n}|w(x,y)|^2(n^2)dydx\bigg)^{1/2} \\
&=n\bigg(\int_{(i-1)/n}^{i/n}\int_{(j-1)/n}^{j/n}|w(x,y)|^2dydx\bigg)^{1/2} \\
\end{align*}
so
\begin{align*}
(CE_1(H^{(n)}))(x)&\le n\sum_{j=1}^n\int_{(i-1)/n}^{i/n}\int_{(j-1)/n}^{j/n}|w(x,y)|^2dydx \\
&=n\int_{(i-1)/n}^{i/n}\int_0^1|w(x,y)|^2dydx \\
&=n\int_{(i-1)/n}^{i/n}(CE_1(H))(x)dx
\le\|CE_1(H)\|_\infty
\end{align*}
and $\|CE_1(H^{(n)})\|_\infty\le\|CE_1(H)\|_\infty$.
Similarly, we get $\|CE_2(H^{(n)})\|_\infty\le\|CE_2(H)\|_\infty$, and~(iii) holds.
\end{proof}

\section{Some quasinilpotent $L^\infty([0,1])$--circular operators}
\label{sec:quasinil}

As mentioned in the introduction, the existence of nontrivial hyperinvariant
subspaces is presently of special interest for quasinilpotent operators in II$_1$--factors.
In this section we give sufficient conditions for an $L^\infty([0,1])$--circular
operator to be quasinilpotent.

\begin{lemma}\label{lem:w}
Let $z_\eta$ be an $L^\infty([0,1])$--circular operator as in Definition~\ref{def:zeta}
and suppose $\eta$ is supported on the set
\[
\{(s,t)\mid0\le s\le t\le1\}
\]
Let $\eps\in(0,1)$, let $\eta_\eps$ be the restriction of $\eta$ to 
\[
\{(s,t)\mid0\le s\le t\le1,\,t-s\le\eps\}
\]
and let $z_{\eta_\eps}$ be the corresponding $L^\infty([0,1])$--circular operator.
Then the spectral radius $r(z_\eta)$ of $z_\eta$ is bounded above by $\|z_{\eta_\eps}\|$.
\end{lemma}
\begin{proof}
Let $\eta_\eps'=\eta-\eta_\eps$.
Then by Proposition~\ref{prop:properties}(ii), $z_\eta$ has the same $*$--moments
as $w+w'$, where $w$ and $w'$ are $L^\infty([0,1])$--circular elements with covariances
$(\alpha_{\eta_\eps},\beta_{\eta_\eps})$ and $(\alpha_{\eta_\eps'},\beta_{\eta_\eps'})$,
respetively, and where $w$ and $w'$ are $*$--free over $L^\infty([0,1])$.
Given $r\in[0,1]$, since
\[
\eta_\eps([r,1]\times[0,r])=0=\eta_\eps'([r-\eps,1]\times[0,r]),
\]
from Lemma~\ref{lem:chiAB}, we have
\[
w\chi_{[0,r]}=\chi_{[0,r]}w\chi_{[0,r]},\qquad w'\chi_{[0,r]}=\chi_{[0,r-\eps]}w'\chi_{[0,r]},
\]
where of course $\chi_{[0,r-\eps]}=0$ if $r\le\eps$.
Letting $p$ be the least integer such that $p\eps\ge1$, we therefore have $(w')^p=0$ and,
given integers $k(0),k(1),\ldots,k(p)\ge0$, we also have
\begin{equation}\label{eq:w'}
w^{k(0)}w'w^{k(1)}w'\cdots w^{k(p-1)}w'w^{k(p)}=0.
\end{equation}
Since $\|z_\eta^n\|=\|(w+w')^n\|$,
by distributing $(w+w')^n$ and using~\eqref{eq:w'}, for $n\ge p$ we obtain
\[
\|z_\eta^n\|\le\sum_{n=0}^{p-1}\binom nq\|w\|^{n-q}\|w'\|^q\le pn^p\max(\|w\|^n,\|w\|^{n-p-1})\max(1,\|w'\|^{p-1}).
\]
Therefore, the spectral radius $r(z_\eta)=\lim_{n\to\infty}\|z_\eta^n\|^{1/n}$,
is bounded above by $\|w\|=\|z_{\eta_\eps}\|$.
\end{proof}

\begin{prop}
Let $z_\eta$ be an $L^\infty([0,1])$--circular operator as in Definition~\ref{def:zeta}.
Suppose $\eta$ is supported on the set
\[
\{(s,t)\mid0\le s\le t\le1\}
\]
and for some $\delta>0$, the restriction of $\eta$ to
\begin{equation}\label{eq:steps}
\{(s,t)\mid0\le s\le t\le1,\,t-s\le\delta\}
\end{equation}
is absolutely continuous with respect to Lebesgue measure and has
bounded Radon--Nikodym derivative.
Then $z_\eta$ is quasinilpotent.
\end{prop}
\begin{proof}
For $0<\eps\le\delta$,
let $\eta_\eps$ be as in Lemma~\ref{lem:w} and
let $H_\eps$ be the Radon-Nikodym derivative of $\eta_\eps$
with respect to Lebesgue measure on $[0,1]^2$.
In this context, equations~\eqref{eq:h1h2} and~\eqref{eq:CEHnorm}
become
\begin{equation}\label{eq:zetaeps}
\|z_{\eta_\eps}\|\le2\max(\|CE_1(H_\eps)\|_\infty,\|CE_2(H_\eps)\|_\infty)^{1/2},
\end{equation}
where
\begin{align*}
CE_1(H_\eps)&=\int_x^{\min(x+\eps,1)}H_\eps(x,y)dy \\
CE_2(H_\eps)&=\int_{\max(0,y-\eps)}^yH_\eps(x,y)dx.
\end{align*}
We thus obtain
\[
\|CE_1(H_\eps)\|_\infty,\|CE_2(H_\eps)\|_\infty\le\eps\|H_\eps\|_\infty\le\eps\|H_\delta\|_\infty.
\]
Consequently, from~\eqref{eq:zetaeps},
$\|z_{\eta_\eps}\|\le2\|H_\delta\|_\infty^{1/2}\sqrt\eps$.
Letting $\eps\to0$ and applying Lemma~\ref{lem:w} yields $r(z_\eta)=0$.
\end{proof}

\bibliographystyle{plain}

\begin{thebibliography}{19}

\bibitem{Br} L.G.\ Brown,
{\em Lidskii's Theorem in the type II case,}
Geometric methods in operator algebras (Kyoto, 1983),
H.\ Araki and E.G.\ Effros, (Eds.),
Pitman Res. Notes Math. Ser.\ {\bf 123}, Longman Sci.\ Tech., 1986, pp.\ 1--35.

\bibitem{Dix} J.\ Dixmier,
{\em Les alg\`ebres d'operateurs dans l'espace hilbertien,}
Gauthier--Villars, 1969.

\bibitem{DP} R.\ Douglas and C.\ Pearcy,
{\em Hyperinvariant subspaces and transitive algebras,}
Mich.\ Math.\ J.\ {\bf 19} (1972), 1-12.

\bibitem{D} K.\ Dykema,
{\em Free products of hyperfinite von Neumann algebras and free dimension,}
Duke Math.\ J.\ {\bf 69} (1993), 97-119.

\bibitem{DH2} K.\ Dykema and U.\ Haagerup,
\emph{DT--operators and decomposability of Voiculescu's circular operator,}
Amer.\ J.\ Math.\ {\bf 126} (2004), 121-189.

\bibitem{DH3} K.\ Dykema and U.\ Haagerup,
\emph{Invariant subspaces of the quasinilpotent DT--operator,}
J.\ Funct.\ Anal.\ {\bf 209} (2004), 332-366.

\bibitem{DS} K.~Dykema and D.~Shlyakhtenko,
{\em Exactness of Cuntz--Pimsner C$^*$--algebras,}
Proc. Edinburgh Math. Soc. {\bf 44} (2001), 425-444. 

\bibitem{F} T.\ Fack,
{\em Sur la notion de valeur charact\'eristique,} J.\ Operator
Theory {\bf 7} (1982), 307-333.

\bibitem{FK} T.\ Fack and H.\ Kosaki,
{\em Generalized $s$--numbers of $\tau$--measurable operators,}
Pacific J.\ Math.\ {\bf 123} (1986), 269-300.

\bibitem{FJKP} C.\ Foia\c s, I.B.\ Jung, E.\ Ko and C.\ Pearcy,
{\em On quasinilpotent operators, III},
preprint.

\bibitem{Ha} U.\ Haagerup,
{\em Spectral decomposition of all operators in a II$_1$--factor which is embeddable in $R^\omega$},
hand--written notes, 2001.

\bibitem{H} T.\ Hoover,
{\em Hyperinvariant subspaces for $n$-normal operators,}
Acta Sci.\ Math.\ (Szeged) {\bf 32} (1971), 109--119.

\bibitem{NSS} A.\ Nica, D.\ Shlyakhtenko and R.\ Speicher,
{\em A characterization of freeness by a factorization property of $R$--transform,}
preprint.

\bibitem{Sh96} D.\ Shlyakhtenko,
{\em Random Gaussian band matrices and freeness with amalgamation,}
Internat.\ Math.\ Res.\ Notices {\bf 1996}, 1013-1025.

\bibitem{Sh97} D.\ Shlyakhtenko,
{\em Free quasi--free states,}
Pacific J.\ Math.\ {\bf 177} (1997), 329-368.

\bibitem{Sh98} D.~Shlyakhtenko,
{\em Gaussian random band matrices and operator--valued free probability theory,}
Quantum probability (Gda\'nsk, 1997), Banach Center Publ., {\bf 43}, Polish Acad.\ Sci., Warsaw, 1998, pp.\ 359--368.

\bibitem{Sh} D.\ Shlyakhtenko,
{\em $A$--valued semicircular systems,}
J.\ Func.\ Anal.\ {\bf 166} (1999), 1-47.

\bibitem{Sn} P.\ \'Sniady,
\emph{Multinomial identities arising from free probability,}
J.\ Combin. Theory A {\bf 101} (2003), 1-19.

\bibitem{Sp} R.\ Speicher,
{\em Combinatorial theory of the free product with amalgamation and operator-valued free probability theory,}
Mem.\ Amer.\ Math.\ Soc.\ {\bf 132} (1998), no.\ 627, x+88 pp.

\bibitem{Tom} J.\ Tomiyama,
{\em On the projection of norm one in W$^*$--algebras,}
Proc.\ Japan Acad.\ {\bf 33} (1957), 608--612.

\bibitem{V} D.\ Voiculescu,
{\em Operations on certain non-commutative operator--valued random variables,}
Recent advances in operator algebras (Orl\'eans, 1992),  Ast\'erisque  No. 232 (1995), pp.\ 243--275.

\bibitem{VDN} D.V.\ Voiculescu, K.J.~Dykema, A.~Nica,
{\em Free Random Variables},
CRM Monograph Series {\bf 1}, American Mathematical Society, 1992.

\end{thebibliography}

\end{document}